\documentclass[conference,onecolumn,12pt]{IEEEtran}
\usepackage[utf8]{inputenc}
\usepackage[a4paper, total={7in, 9in}]{geometry}
\usepackage{listings}
\usepackage{parskip}
\usepackage{hyperref}
\newtheorem{lemma}{Lemma}
\newtheorem{theorem}{Theorem}
\newtheorem{remark}{Remark}
\usepackage{algorithm}
\usepackage{algorithmic}
\usepackage{multirow}
\usepackage{cite}
\usepackage[document]{ragged2e}
\usepackage{color} %red, green, blue, yellow, cyan, magenta, black, white
\definecolor{mygreen}{RGB}{28,172,0} % color values Red, Green, Blue
\definecolor{mylilas}{RGB}{170,55,241}
\usepackage{graphicx}
\usepackage{amsmath, amssymb, txfonts, amsfonts}
\AtBeginDocument{\pagenumbering{arabic}}

\date{}
\title{\textbf{Robust Maximum Correntropy Kalman Filter}}
\def\BibTeX{{\rm B\kern-.05em{\sc i\kern-.025em b}\kern-.08em
		T\kern-.1667em\lower.7ex\hbox{E}\kern-.125emX}}	
\author{Joydeb Saha \textsuperscript{1},  Shovan Bhaumik\textsuperscript{2},\\ 
	\{\IEEEauthorblockN{\textsuperscript{1}joydeb$ \_ $2121ee32, ~ \textsuperscript{2}shovan.bhaumik\} @iitp.ac.in \\\textsuperscript{1,2} Indian Institute of Technology Patna (IITP),\\ Patna, Bihar, India}
	\and
	\\ \IEEEauthorblockN{}
}

\begin{document}

\maketitle 
\justifying

\begin{abstract}             
	The Kalman filter provides an optimal estimation for a linear system with Gaussian noise. However when the noises are non-Gaussian in nature, its performance deteriorates rapidly. For non-Gaussian noises, maximum correntropy Kalman filter (MCKF) is developed which provides an improved result. But when the system model differs from nominal consideration, the performance of the MCKF degrades. For such cases, we have proposed a new robust filtering technique which maximize a cost function defined by exponential of weighted past and present errors along with the Gaussian kernel function. By solving this cost criteria we have developed prior and posterior mean and covariance matrix propagation equations. By maximizing the correntropy function of error matrix, we have selected the kernel bandwidth value at each time step. Further the conditions for convergence of the proposed algorithm is also derived. Two numerical examples are presented to show the usefulness of the new filtering technique.
\end{abstract}

%\begin{keyword}
%Maximum correntropy kalman filter (MCKF), Robust maximum correntropy kalman filter (RMCKF), Maximum correntropy criteria (MCC), System uncertainty model 
%\end{keyword}

%===============================================================================

\section{Introduction}
State estimation is a very important technique used in various industrial problems and in research applications such as target tracking, navigation system, communication system, image processing, system identification, data fusion, satellite state estimation, and many more. The Kalman filter provides and optimal estimation of states for linear systems when the noises are Gaussian in nature. But for non-Gaussian noises, the performance of the Kalman filter degrades drastically. \\
To resolve this limitation, a few approaches such as minimum error entropy based Kalman filter \cite{chen2019minimum}, Bayesian inference algorithm \cite{zhang2018multipath}, maximum correntropy Kalman filter (MCKF) \cite{chen2017maximum} \emph{etc.} are developed. Similarity is a key concept to express or measure the quantity of a temporal signal. Correntropy is directly related to the probability of the similarity of two random variables in a neighborhood of a joint space defined by kernel bandwidth. It is a bivariate function that produces a scalar which contains second and higher order pdf moments. In recent years, correntropy \cite{liu2007correntropy} based filtering is being used for state estimation in presence of non-Gaussian noises where correntropy is maximized and such filters are known maximum correntropy Kalman filter (MCKF) \cite{chen2017maximum}. In traditional Kalman filter, mean square error is minimized which deals with second order pdf moment whereas maximum correntropy criteria (MCC) considers all the higher order even moments along with it. This is the main reason why MCC based filters give better results in presence of non-Gaussian noise than traditional Kalman filter. \\ 
Kernel function plays an important role in correntropy based filtering techniques. Gaussian kernel is very popularly used in MCKF.  In literature, few more kernels such as Laplacian kernel \cite{hu2021robust}, Gaussian mixture kernel \cite{wang2022mixture} are also available. But Gaussian kernel is smooth, symmetric and integral of product of two Gaussians remains Gaussian \cite{liu2007correntropy}. Because of these properties, Gaussian kernels are preferred. \\
In many practical applications, we do not know the process model with certainty. In such a case, performance of MCKF may degrade for model mismatch. We need a robust algorithm to handle this scenario. In literature, risk sensitive filters (RSF) are present with Gaussian noise consideration. In \cite{bhaumik2005risk}, \cite{bhaumik2009risk} and \cite{bhaumik2008improved}, a detailed study on the formulation of risk sensitive estimation problem is explained. But for non-Gaussian noise with system uncertainty, nothing is available in literature. By merging the concept of risk sensitive filter with the idea of correntropy filter, we have formulated a new cost function which is based on weighted sum of all the past errors and weighted present error having the Gaussian kernel function with kernel bandwidth. By maximizing this, a new algorithm is formulated that is a good fit for system uncertainty model with non-Gaussian noise. \\
In correntropy based filters, kernel bandwidth owns significant importance in the performance of the filtering technique. Selection of the proper bandwidth value is a major challenge to the researchers. Few publications discussed regarding the adaptive kernel bandwidth selection approach \cite{cinar2012hidden, hou2017maximum, fakoorian2019maximum} but these do not guarantee the optimal value. We have proposed an alternative cost function using Gaussian kernel to numerically select the bandwidth value for each time step. We also derived the convergence and stability criteria for our proposed filter. 

\section{Problem Formulation}
Let us consider a linear system having the following process and measurement equations:
\begin{equation}
	\mathcal{X}_{k+1} = (F_{k} + \Delta F_{k}) \mathcal{X}_{k} + q_{k}, \label{process_equ}
\end{equation}
\begin{equation}
	\mathcal{Y}_{k} = H_{k} \mathcal{X}_{k} + r_{k}, \label{measure_equ}
\end{equation}

where $\mathcal{X}_{k} = \begin{bmatrix}
	x_{1,k} & x_{2,k} & . & . & . & x_{n,k}
\end{bmatrix}^T$ is the state vector of the system, $\mathcal{Y}_{k} = \begin{bmatrix}
	y_{1,k} & y_{2,k} & . & . & . & y_{m,k}
\end{bmatrix}^T$ is the measurement vector. $F_{k}$ and $H_{k}$ are the state and measurement matrices respectively. $\Delta F_{k}$ is an arbitrary and deterministic unknown parameter in process model which defines the uncertainty of the system. We assume $\Delta F_{k}$ is bounded in such a way so that the perturbed system remains stable. If the system parameters are accurately known, $\Delta F_{k}=0$. Process noise $q_{k}$ and measurement noise $r_{k}$ are zero mean and follow a non-Gaussian distribution with equivalent covariance $Q_{k}$ and $R_{k}$ respectively \emph{i.e.} $q_{k} = \sum a_{i}\mathcal{N} (0, Q_{i})$ and $r_{k} = \sum a_{i}\mathcal{N} (0, R_{i})$. We also consider the noises are uncorrelated to each other \emph{i.e.} $\mathbf{E}[q_{k} {r_{k}}^T] = 0$. Our objective is to find the posterior estimate $\mathcal{\hat{X}}_{k|k}$ from the measurements $\mathcal{Y}_{1:k}$ \emph{i.e.} $p(\mathcal{\hat{X}}_{k|k}|\mathcal{Y}_{1:k})$ for the system defined in (\ref{process_equ}) and (\ref{measure_equ}) where filter assumes the system without perturbation \emph{i.e.} $\Delta F_{k} = 0$.

\begin{remark}
	The noises here are $q_{k}$ and $r_{k}$ which are non-Gaussian. They can be expressed as a weighted sum of many Gaussian noises in the form of $\sum a_{i} \mathcal{N} (\nu_{i}, \Sigma_{i})$, where $a_{i}$ are weights satisfying the condition $\sum a_{i} = 1$. $\nu_{i}$ and $\Sigma_{i}$ are respectively mean and covariance matrics of the $i^{th}$ normal distributions.
\end{remark}

\section{Correntropy} \label{sec:correntropy}
Correntropy is directly related to the probability of the similarity of two random variables in a neighbourhood of the joint space defined by kernel bandwidth. Correntropy function produces a scalar which contains second and higher order pdf moments. We denote $p(x_{i}, \hat{x}_{i})$ and $F(x_{i}, \hat{x}_{i})$ as joint pdf and CDF respectively of each $i^{th}$ state, $i \in [1,n]$. So, the correntropy of each state $x_{i}$ can be defined as
\begin{equation} 
	\begin{split}
		\mathsf{v}(x_i, \hat{x}_i) = &~ \mathbf{E}[\mathsf{k}(x_i, \hat{x}_i)] \\
		= &~ \int \int \mathsf{k}(x_i, \hat{x}_i) p(x_i, \hat{x}_i) dx_i d\hat{x}_i \\
		= &~ \int \mathsf{k}(x_i, \hat{x}_i) dF(x_i, \hat{x}_i), \label{corr_func}
	\end{split}
\end{equation}
where $\mathbf{E}$ is the expectation operator and $\mathsf{k}(.,.)$ is the kernel function.
For Gaussian kernel, $\mathsf{k}(x_i, \hat{x}_i)$ can be expressed as 
\begin{equation}
	\mathsf{k}(x_i, \hat{x}_i) =  G_{\sigma}(x_i, \hat{x}_i)
	=  \text{exp}(-\frac{(x_i - \hat{x}_i)^2}{2 \sigma^2}), \label{gau_ker_func}
\end{equation}
where $\sigma$ defines the kernel bandwidth. Now, considering the error function, the Gaussian kernel can be written as 
\begin{equation}	
	G_{\sigma}(x_i, \hat{x}_{i}) = G_{\sigma}(e_{i}) 
	= \text{exp}(-\frac{e_{i}^2}{2\sigma^2}). \label{gau_ker}
\end{equation}
Hence, for Gaussian kernel the correntropy of each state of the system will be 
\begin{equation}
	\mathsf{v}_{\sigma}(x_i, \hat{x}_i) = \int	G_{\sigma}(e_{i}) dF(x_i, \hat{x}_i). \label{gau_corr_fun}
\end{equation}
In practical cases, availability of sample data are limited for which joint CDF $F(x_i, \hat{x_i})$ is usually unavailable. So, we can write the correntropy for each state with the help of sample mean estimator. So the total correntropy of the system at any time step $k$ will be
\begin{equation}
	\hat{\mathsf{v}}_{\sigma}(\mathcal{X}_{k}, \mathcal{\hat{X}}_{k}) = \frac{1}{n} \sum_{i=1}^{n} G_{\sigma}(e_{k,i}). \label{gau_corr_time_step}
\end{equation}
where $e_{k,i} = (x_{k,i} - \hat{x}_{k,i})$.

\begin{remark}
	A few important properties of correntropy function can be explained if we expand (\ref{gau_corr_time_step}) by the Taylor series. With the Taylor series expansion, we get $\hat{\mathsf{v}}_{\sigma}(x_i, \hat{x}_i) = \sum_{N=0}^{\infty} \frac{(-1)^N}{2^N \sigma^{2N} N!} \mathbf{E}[(x_i - \hat{x}_i)^{2N}]$. it can be said that for the Gaussian kernel function, correntropy becomes the sum of all even moments of the difference between two random variables. This creates a major difference between the mean square error (MSE) and the correntropy criteria as MSE deals with second order moment only.
\end{remark}
\begin{remark}
	For Gaussian distribution, MSE provides the optimal estimation as second order moment is sufficient for that case. But to describe a non-Gaussian distribution, only second order moment is not enough. That's why MSE based filtering approach such as Kalman filter fails to provide an optimal solution in case of non-Gaussian distribution and we look towards correntropy based filtering approach.
\end{remark}

\begin{remark}
	The kernel bandwidth $\sigma$ works as a weighting parameter to the even order monents. Hence, it can be said that increment in $\sigma$ will impact more in higher order moments as compared to second order moment. For a large value of kernel bandwidth, higher order moments will be near to zero. Hence it will work like MSE criteria.
\end{remark}

\section{Robust Maximum Correntropy Kalman Filter (RMCKF)}
\subsection{Cost Function} \label{sec:cf}
In this section, we define a cost function for robust MCKF which is different than existing cost functions.
To describe cost function, first we augment the system (\ref{process_equ}) and  (\ref{measure_equ}) as follows,
\begin{equation}
	\begin{bmatrix}
		\mathcal{\hat{X}}_{k|k-1} \\
		\mathcal{Y}_{k}
	\end{bmatrix} = 
	\begin{bmatrix}
		I \\
		H_{k}
	\end{bmatrix} \mathcal{X}_{k} + v_{k}, \label{modf_sys_mod}
\end{equation}
where \begin{equation}
	v_{k} = \begin{bmatrix}
		-(\mathcal{X}_{k} - \mathcal{\hat{X}}_{k|k-1}) \\
		r_{k}
	\end{bmatrix}.
\end{equation}
So, \begin{equation}
	\begin{split}
		\mathbf{E}[v_{k} v_{k}^T] = &~ \begin{bmatrix}
			P_{k|k-1} & 0 \\
			0 & R_{k}
		\end{bmatrix} \\
		= &~ \begin{bmatrix}
			B_{p, k|k-1} B_{p, k|k-1}^T & 0 \\
			0 & B_{r,k} B_{r,k}^T
		\end{bmatrix} \\
		= &~ B_{k} B_{k}^T.
	\end{split} \label{v_k_cov}
\end{equation}
$B_{P, k|k-1}$ and $B_{r,k}$ are square-roots of $P_{k|k-1}$ and $R_{k}$ respectively. From (\ref{v_k_cov}) it can be said that $B_{k}$ is the square-root of the matrix $\begin{bmatrix}
	P_{k|k-1} & 0 \\
	0 & R_{k}
\end{bmatrix}$
and  can be expressed as a diagonal matrix
\begin{equation}
	B_{k} = \begin{bmatrix}
		B_{p,k|k-1} & 0 \\
		0 & B_{r,k}
	\end{bmatrix}.
\end{equation}
Now, left multiplying both sides of (\ref{modf_sys_mod}) by $B_{k}^{-1}$, we will get
\begin{equation}
	D_{k} = W_{k} \mathcal{X}_{k} + e_{k}, \label{err_mat_def}
\end{equation}
where $D_{k} = B_{k}^{-1} \begin{bmatrix}
	\mathcal{\hat{X}}_{k|k-1} &
	\mathcal{Y}_{k}
\end{bmatrix}^T$,
$W_{k} = B_{k}^{-1} \begin{bmatrix}
	I &
	H_{k}
\end{bmatrix}^T$ and $e_{k} = B_{k}^{-1} v_{k}$. $e_{k}$ defines the error matrix of dimension $(n+m) \times 1$. \\
Now, we define a cost function as follows,
\begin{equation}
	J_{L}(\mathcal{X}_{k}) = \frac{1}{L} \sum_{i=1}^{L} [\text{exp}\{- \sum_{j=1}^{k-1} \frac{\mu_{1} e_{j,i}^2}{2 \sigma^2} -   \frac{\mu_{2} e_{k,i}^2}{2 \sigma^2}\}],\label{RSMCC_cost_criteria}
\end{equation}
where $e_{j,i}$ is the $i^{th}$ element of the error matrix at $j^{th}$ time step $(e_{j})$. $\mu_{1}$ and $\mu_{2}$ are two risk sensitive parameters, $\sigma$ is the kernel bandwidth, and $L = n+m$. Using (\ref{err_mat_def}), $e_{k,i}$ can be defined as
\begin{equation} \label{err_ele}
	e_{k,i} = d_{k,i} - w_{k,i} x_{k,i},
\end{equation}
where $e_{k,i}$, $d_{k,i}$ and $x_{k,i}$ are $i^{th}$ element of $e_{k}$, $D_{k}$ and $\mathcal{X}_{k}$ respectively and $w_{k,i}$ is $i^{th}$ row of $W_{k}$. Our objective is to find an optimal posterior estimate of state $\mathcal{\hat{X}}^{*}_{k|k}$ from the received measurements $\mathcal{Y}_{1:k}$ by maximizing the cost function (\ref{RSMCC_cost_criteria}) that is
\begin{equation}
	\mathcal{\hat{X}}^{*}_{k|k} = \text{arg} \text{ $\displaystyle \max_{\mathcal{X}_{k}}$} J_{L}(\mathcal{X}_k). \label{arg_max_equ}
\end{equation}

\begin{remark}
	In (\ref{RSMCC_cost_criteria}), we propose a new cost function which does not exist in earlier literature. It can be thought of a combination of maximum correntropy criteria \cite{chen2017maximum} and risk sensitive cost function \cite{bhaumik2005risk} Please note that the term $[\text{exp}\{- \sum_{j=1}^{k-1} \frac{\mu_{1} e_{j,i}^2}{2 \sigma^2} - \frac{\mu_{2} e_{k,i}^2}{2 \sigma^2}\}]$ is a modified form of Gaussian kernel $G_{\sigma}(e)$ described in (\ref{gau_ker}) where the cost function $J_{L}(\mathcal{X}_{k})$ is the correntropy function.
\end{remark}

\begin{remark}
	One notable point is that the error matrix $e_{k}$ is not the direct difference between the true state and estimated state. Rather there is weighting factors $B_{p,k|k-1}^{-1}$ and $B_{r,k}^{-1}$ for process error and for measurement error respectively. This is a fundamental difference in the construction of error matrix in correntropy based filter and in MSE based filter. 
\end{remark}

\begin{remark}
	For the risk parameters $\mu_{1} = 0$, $\mu_{2} = 1$ and finite kernel bandwidth \emph{i.e.} $\sigma \neq \infty$, the cost function becomes the same as maximum correntropy cost function as mentioned in \cite{chen2017maximum}. For $\mu_{1} \neq 0$ and $\mu_{2} = 1$ and infinite kernel bandwidth \emph{i.e.} $\sigma \to \infty$, the cost function becomes risk sensitive cost function as defined in section 2 of \cite{bhaumik2005risk}. For $\mu_{1} = 0$, $\mu_{2} = 1$ and $\sigma \to \infty$, the cost function is same as mean square error cost function which leads to the KF when we minimize it.
\end{remark}

\subsection{Formulation of Robust Maximum Correntropy Kalman Filter (RMCKF)}
The posterior information state density, $\xi_{k}$ is defined as $\xi_{k} = p(\mathcal{X}_{k}| \mathcal{I}_{k})$ where $\mathcal{I}_{k} = \{\mathcal{Y}_{1:k}, e_{1|1}, \cdots , e_{k-1|k-1}\}$. Further, using Eqn. (17) of \cite{tiwari2022risk}, the information state density is further written as 
\begin{equation}
	\xi_{k} = p(\mathcal{X}_{k}| \mathcal{I}_{k}) = \text{exp}(- \sum_{j=1}^{k-1} \frac{\mu_{1} e_{j,i}^2}{2 \sigma^2}) p (\mathcal{X}_{k}|\mathcal{Y}_{1:k}). \label{information_state}
\end{equation}

\begin{lemma} 
	The expression of recursive update of $\xi_{k}$ is 
	\begin{equation}
		\xi_{k} = \zeta_{k} p(\mathcal{Y}_{k}|\mathcal{X}_{k}) \int p(\mathcal{X}_{k}|\mathcal{X}_{k-1}) \text{exp}(-  \frac{\mu_{1} e_{k-1,i}^2}{2 \sigma^2}) \xi_{k-1} d\mathcal{X}_{k-1}.
	\end{equation}
	where $\zeta_{k} = \frac{1}{p(\mathcal{Y}_{k}|\mathcal{Y}_{1:k-1})}$.
\end{lemma}
%\begin{pf}
\textbf{Proof:}
Using the Bayes' theorem, the posterior probability density function (pdf) of states, $p(\mathcal{X}_{k}|\mathcal{Y}_{1:k})$ can be written as
\begin{equation}
	\begin{split}
		p(\mathcal{X}_{k}|\mathcal{Y}_{1:k}) = &~ \frac{p(\mathcal{Y}_{1:k}|\mathcal{X}_{k}) p(\mathcal{X}_{k})}{p(\mathcal{Y}_{1:k})} \\
		= &~ \frac{p(\mathcal{Y}_{k}|\mathcal{Y}_{1:k-1},\mathcal{X}_{k}) p(\mathcal{Y}_{1:k-1}|\mathcal{X}_{k}) p(\mathcal{X}_{k})}{p(\mathcal{Y}_{k}|\mathcal{Y}_{1:k-1}) p(\mathcal{Y}_{1:k-1})} \\
		= &~ \frac{p(\mathcal{Y}_{k}|\mathcal{Y}_{1:k-1},\mathcal{X}_{k}) p(\mathcal{X}_{k}|\mathcal{Y}_{1:k-1}) p(\mathcal{Y}_{1:k-1}) p(\mathcal{X}_{k})}{p(\mathcal{Y}_{k}|\mathcal{Y}_{1:k-1}) p(\mathcal{Y}_{1:k-1})p(\mathcal{X}_{k})} \\
		= &~ \frac{p(\mathcal{Y}_{k}|\mathcal{X}_{k}, \mathcal{Y}_{1:k-1}) p(\mathcal{X}_{k}|\mathcal{Y}_{1:k-1}) }{p(\mathcal{Y}_{k}|\mathcal{Y}_{1:k-1})} \\
		= &~ \zeta_{k} p(\mathcal{Y}_{k}|\mathcal{X}_{k}, \mathcal{Y}_{1:k-1}) p(\mathcal{X}_{k}|\mathcal{Y}_{1:k-1}).
	\end{split} \label{post_den}
\end{equation}
Here, $\zeta_{k} = \frac{1}{p(\mathcal{Y}_{k}|\mathcal{Y}_{1:k-1})}$ is a normalizing constant, where $p(\mathcal{Y}_{k}|\mathcal{Y}_{1:k-1}) = \int p (\mathcal{Y}_{k}|\mathcal{X}_{k}) p (\mathcal{X}_{k}|\mathcal{Y}_{1:k-1}) d\mathcal{X}_{k}$. Considering $\mathcal{Y}_{k}$ is independent from the previous measurements $\mathcal{Y}_{1:k-1}$, we can write $p(\mathcal{Y}_{k}|\mathcal{X}_{k}, \mathcal{Y}_{1:k-1}) = p(\mathcal{Y}_{k}|\mathcal{X}_{k})$. Hence, Eqn.(\ref{post_den}) becomes
\begin{equation}
	p(\mathcal{X}_{k}|\mathcal{Y}_{1:k}) = \zeta_{k} p(\mathcal{Y}_{k}|\mathcal{X}_{k}) p(\mathcal{X}_{k}|\mathcal{Y}_{1:k-1}). \label{post_den_upd}
\end{equation}
Applying Chapman-Kolmogorov integral, $p(\mathcal{X}_{k}|\mathcal{Y}_{1:k})$ can be expressed as 
\begin{equation}
	p(\mathcal{X}_{k}|\mathcal{Y}_{1:k}) = \zeta_{k} p(\mathcal{Y}_{k}|\mathcal{X}_{k}) \int p(\mathcal{X}_{k}|\mathcal{X}_{k-1}) p (\mathcal{X}_{k-1}|\mathcal{Y}_{1:k-1}) d\mathcal{X}_{k-1}. \label{post_den_state}
\end{equation}

Now, substituting $p (\mathcal{X}_{k}|\mathcal{Y}_{1:k})$ in (\ref{information_state}) from (\ref{post_den_state}) we get
\begin{equation}
	\begin{split}
		\xi_{k} = p(\mathcal{X}_{k}| \mathcal{I}_{k}) &~
		= \text{exp}(-\sum_{j=1}^{k-1} \frac{\mu_{1} e_{j,i}^2}{2 \sigma^2}) \times \zeta_{k} p(\mathcal{Y}_{k}|\mathcal{X}_{k}) \\
		&~ \times \int p(\mathcal{X}_{k}|\mathcal{X}_{k-1}) p (\mathcal{X}_{k-1}|\mathcal{Y}_{1:k-1}) d\mathcal{X}_{k-1} \\
		= &~ \zeta_{k} p(\mathcal{Y}_{k}|\mathcal{X}_{k}) \int p(\mathcal{X}_{k}|\mathcal{X}_{k-1}) \text{exp}(-  \frac{\mu_{1} e_{k-1,i}^2}{2 \sigma^2}) \\
		&~ \times \text{exp}(- \sum_{j=1}^{k-2} \frac{\mu_{1} e_{j,i}^2}{2 \sigma^2}) p (\mathcal{X}_{k-1}|\mathcal{Y}_{1:k-1}) d\mathcal{X}_{k-1} \\
		= &~ \zeta_{k} p(\mathcal{Y}_{k}|\mathcal{X}_{k}) \int p(\mathcal{X}_{k}|\mathcal{X}_{k-1}) \text{exp}(-  \frac{\mu_{1} e_{k-1,i}^2}{2 \sigma^2}) \xi_{k-1} d\mathcal{X}_{k-1}.
	\end{split} \label{inter_inf_st_n}
\end{equation}

\begin{flushright}
	$\blacksquare$	
\end{flushright}
%\end{pf}

\begin{remark}
	The cost function described in (\ref{RSMCC_cost_criteria}) can alternatively be expressed with information state pdf as
	\begin{equation}
		J_{L}(\mathcal{X}_{k}) =\int \text{exp}(-\sum_{j=1}^{k-1} \frac{\mu_{1} e_{j,i}^2}{2 \sigma^2}) \times \text{exp} (-\frac{\mu_{2} e_{k,i}^2}{2\sigma^2}) \times p(\mathcal{X}_{k}|\mathcal{Y}_{1:k}) d\mathcal{X}_{k}, \label{inter_cost}
	\end{equation}
	or,
	\begin{equation}
		J_{L}(\mathcal{X}_{k}) = \int \text{exp}(- \sum_{j=1}^{k-1} \frac{\mu_{1} e_{j,i}^2}{2 \sigma^2}) \xi_{k} d\mathcal{X}_{k}.
	\end{equation}
\end{remark}

\begin{remark} \label{rem:Gauassm}
	Following \cite{tiwari2022risk}, we define $\int p(\mathcal{X}_{k}|\mathcal{X}_{k-1}) \text{exp}(-  \frac{\mu_{1} e_{k-1,i}^2}{2 \sigma^2}) \xi_{k-1} d\mathcal{X}_{k-1}$ as prior information state density and symbolized it as  $p(\mathcal{X}_{k}| \mathcal{I}_{k-1}, e_{k-1|k-1})$. From the Lemma 1, we see that the information state pdf does not remain Gaussian even if we begin with a Gaussian information state. Because the likelihood and state transition density become non-Gaussian due to the presence of non-Gaussian process and measurement noises. However, here we approximated $\xi_{k}$ as Gaussian with a mean and equivalent covariance.
\end{remark}

% The initial point of the recursion is to be considered as $\xi_{0} = p(\mathcal{X}_{0})$. So, (\ref{inter_inf_st_n}) can be written  in the form of (\ref{xi_k}).

\begin{theorem}
	Under the assumption of remark \ref{rem:Gauassm}, the expressions of prior mean and prior error covariance are
	\begin{equation}
		\mathcal{\hat{X}}_{k|k-1} = F_{k} \mathcal{\hat{X}}_{k-1|k-1},  \label{Pri_state}
	\end{equation}
	\begin{equation}
		P_{k|k-1} = F_{k-1} (P_{k-1|k-1}^{-1} - 2 \mu_{1} I)^{-1} F_{k-1}^T + Q_{k-1} \label{Pri_cov}.
	\end{equation}
	
\end{theorem}
%\begin{pf}
\textbf{Proof:}
The prior information states which are assumed as Gaussian can be expressed as following:
\begin{equation}
	\xi_{k-1} = \zeta_{k} (2 \pi)^{-n/2} |P_{k-1|k-1}|^{-1/2} \text{exp}(-\frac{1}{2} (\mathcal{X}_{k-1} - \mathcal{\hat{X}}_{k-1|k-1}) P_{k-1|k-1}^{-1} (\mathcal{X}_{k-1} - \mathcal{\hat{X}}_{k-1|k-1})^T).
\end{equation}
Now, substituting the value of $\xi_{k-1}$ in prior information state density, we receive
\begin{equation}
	\begin{split}
		&~ p(\mathcal{X}_{k}| \mathcal{I}_{k-1}, e_{k-1|k-1}) \\
		= &~  \int p(\mathcal{X}_{k}|\mathcal{X}_{k-1}) \text{exp}(-  \frac{\mu_{1} e_{i,k-1}^2}{2 \sigma^2}) \zeta_{k-1} (2 \pi)^{-n/2} |P_{k-1|k-1}|^{-1/2} \\
		&~ \times  \text{exp}(-\frac{1}{2} (\mathcal{X}_{k-1} - \mathcal{\hat{X}}_{k-1|k-1})P_{k-1|k-1}^{-1} (\mathcal{X}_{k-1} - \mathcal{\hat{X}}_{k-1|k-1})^T) d\mathcal{X}_{k-1} \\
		= &~ \zeta_{k-1} (2 \pi)^{-n/2} [|P_{k-1|k-1}|^{-1/2} \text{exp}(\frac{1}{2 \sigma^2}) \int p(\mathcal{X}_{k}|\mathcal{X}_{k-1}) \\
		&~ \times \text{exp}(-\frac{1}{2} (\mathcal{X}_{k-1} - \mathcal{\hat{X}}_{k-1|k-1})P_{k-1|k-1}^{-1} (\mathcal{X}_{k-1} - \mathcal{\hat{X}}_{k-1|k-1})^T) \\
		&~ \times\text{exp}(-\frac{1}{2} (\mathcal{X}_{k-1} - \mathcal{\hat{X}}_{k-1|k-1}) (-2 \mu_{1} I) (\mathcal{X}_{k-1} - \mathcal{\hat{X}}_{k-1|k-1})^T)] d\mathcal{X}_{k-1} \\
		= &~ \zeta_{k-1} (2 \pi)^{-n/2} \text{exp}(\frac{1}{2 \sigma^2}) [|P_{k-1|k-1}|^{-1/2}  \int p(\mathcal{X}_{k}|\mathcal{X}_{k-1}) \\
		&~ \times \text{exp}(-\frac{1}{2} (\mathcal{X}_{k-1} - \mathcal{\hat{X}}_{k-1|k-1})(P_{k-1|k-1}^{-1} - 2 \mu_{1}I) (\mathcal{X}_{k-1} - \mathcal{\hat{X}}_{k-1|k-1})^T)] \\
		&~ d\mathcal{X}_{k-1},
	\end{split} \label{inf_pri}
\end{equation}

where $(P_{k-1|k-1}^{-1} - 2 \mu_{1}I)$ should be invertible. It is obvious that the (\ref{inf_pri}) represents a Gaussian distribution with mean $\mathcal{\hat{X}}_{k-1|k-1}$ and covariance $(P_{k-1|k-1}^{-1} - 2 \mu_{1}I)^{-1}$. Hence, (\ref{Pri_state}) and (\ref{Pri_cov}) are obtained.
\begin{flushright}
	$\blacksquare$	
\end{flushright}
%\end{pf}

\begin{theorem}
	The expression of the posterior estimate of state and posterior error covariance will be
	\begin{equation}
		\mathcal{\hat{X}}_{k|k} = \mathcal{\hat{X}}_{k|k-1} + K_{k} (\mathcal{Y}_{k} - H_{k}\mathcal{\hat{X}}_{k|k-1}), \label{Post_state}
	\end{equation} 
	
	\begin{equation}
		P_{k|k} = (I - K_{k} H_{k}) P_{k|k-1} (I - K_{k} H_{k})^T + K_{k} R_{k} K_{k}^T, \label{post_err_cov}
	\end{equation} 	
	
	where
	\begin{equation}
		K_{k} = \bar{P}_{k|k-1} H_{k}^T (H_{k} \bar{P}_{k|k-1} H_{k}^T + \bar{R}_{k})^{-1}, \label{K_value}
	\end{equation}
	$\bar{P}_{k|k-1} = (B_{p,k|k-1}) \Pi_{p,k}^{-1} (B_{p,k|k-1})^T$ , $\bar{R}_{k} = (B_{r,k}) \Pi_{r,k}^{-1} (B_{r,k})^T$,
	$\Pi_{p,k} = diag 
	(\Pi_{p,k,1}, \; \Pi_{p,k,2}, \cdots , \Pi_{p,k,n})
	$, and $\Pi_{r,k} = diag
	(\Pi_{r,k,1}, \; \Pi_{r,k,2}, \cdots , \Pi_{r,k,m})$ with $\Pi_{p,k,i} = \text{exp}(\rho_{p,i} -\frac{\mu_{2} e_{p,k,i}^2} {2 \sigma^2})$ and $\Pi_{r,k,i} = \text{exp}(\rho_{r,i} -\frac{\mu_{2} e_{r,k,i}^2} {2 \sigma^2})$, $e_{k} = \begin{bmatrix}
		e_{p,k} & e_{r,k}
	\end{bmatrix}^T$ where $e_{p,k} = -B_{p,k|k-1}^{-1} (\mathcal{X}_{k} - \mathcal{\hat{X}}_{k|k-1})$ and $e_{r,k} = B_{r,k}^{-1} (\mathcal{Y}_{k} - H_{k}\mathcal{X}_{k})$ represent weighted process and measurement errors respectively. And $\rho_{p,i} = \sum_{j = 1}^{k-1} (- \frac{\mu_{1} e_{p,j,i}^2}{2 \sigma^2})$ and $\rho_{r,i} = \sum_{j = 1}^{k-1} (- \frac{\mu_{1} e_{r,j,i}^2}{2 \sigma^2})$ denote the weighted past process and measurement errors respectively.
\end{theorem}

%\begin{pf}
\textbf{Proof:}
Partially differentiating $e_{p,k}$ and $e_{r,k}$ \emph{w.r.t.} $\mathcal{X}_{k}$ the following expansions will occur
\begin{equation}
	\frac{\partial e_{p,k}}{\partial \mathcal{X}_{k}} = -(B_{p,k|k-1}^{-1})^T, \label{deriv_ep}
\end{equation} and
\begin{equation}
	\frac{\partial e_{r,k}}{\partial \mathcal{X}_{k}}
	= -H_{k}^T (B_{r,k}^{-1})^T. \label{deriv_er}
\end{equation}

The cost function described in (\ref{RSMCC_cost_criteria}) can further be written as 
\begin{equation}
	\begin{split}
		J_{L}(\mathcal{X}_{k}) &~ = \frac{1}{L} \sum_{i = 1}^{L} \text{exp}(\rho_{i} -\frac{\mu_{2} e_{k,i}^2}{2 \sigma^2}) \\
		&~ = \frac{1}{L} [\sum_{i = 1}^{n} \text{exp}(\rho_{i} -\frac{\mu_{2} e_{p,k,i}^2}{2 \sigma^2}) + \sum_{i = 1}^{m} \text{exp}(\rho_{i} -\frac{\mu_{2} e_{r,k,i}^2}{2 \sigma^2})]
	\end{split}
\end{equation} 
where $\rho_{i} = \sum_{j = 1}^{k-1} (- \frac{\mu_{1} e_{j,i}^2}{2 \sigma^2})$. Now, partially differentiating $J_{L}(\mathcal{X}_{k})$ \emph{w.r.t.} $\mathcal{X}_{k}$, we will get
\begin{equation}
	\begin{split}
		\frac{\partial J_{L}(\mathcal{X}_k)}{\partial \mathcal{X}_{k}} = &~ \frac{1}{L} [\sum_{i=1}^{n} \text{exp}(\rho_{p,i}-\frac{\mu_{2} e_{p,k,i}^2}{2 \sigma^2})(- \frac{\mu_{2} e_{p,k,i}}{\sigma^2}) \frac{\partial e_{p,k,i}}{\partial \mathcal{X}_k}
		+ \sum_{i=1}^{m} \text{exp}(\rho_{r,i}-\frac{\mu_{2} e_{r,k,i}^2}{2 \sigma^2})(- \frac{\mu_{2} e_{r,k,i}}{\sigma^2}) \frac{\partial e_{r,k,i}}{\partial \mathcal{X}_k}] \\
		= &~ 0. \label{diff_J}
	\end{split}
\end{equation}
By simplifying (\ref{diff_J}), we will get

\begin{equation}
	\sum_{i=1}^{n} \text{exp}(\rho_{p,i} - \frac{\mu_{2} e_{p,k,i}^2}{2\sigma^2})e_{p,k,i} \frac{\partial e_{p,k,i}}{\partial \mathcal{X}_k}
	+ \sum_{i=1}^{m} \text{exp}(\rho_{r,i} - \frac{\mu_{2} e_{r,k,i}^2}{2\sigma^2})e_{r,k,i} 
	\frac{\partial e_{r,k,i}}{\partial \mathcal{X}_k}=0. \label{deriv_cost}
\end{equation}

From (\ref{deriv_ep}), (\ref{deriv_er}) and (\ref{deriv_cost}) following recursive equation can be obtained:
\begin{equation}
	(B_{p,k|k-1}^{-1})^T \Pi_{p,k} B_{p,k|k-1}^{-1} (\hat{\mathcal{X}}_{k|k} - \mathcal{\hat{X}}_{k|k-1})
	- H_{k}^T (B_{r,k}^{-1})^T \Pi_{r,k} B_{r,k}^{-1} (\mathcal{Y}_{k} - H_{k}\hat{\mathcal{X}}_{k|k}) = 0, \label{recur_equ}
\end{equation}
Considering  $B_{p,k|k-1} \Pi_{p,k}^{-1} B_{p,k|k-1}^T = \bar{P}_{k|k-1}$ and $B_{r,k} \Pi_{r,k}^{-1} B_{r,k}^T = \bar{R}_{k}$, (\ref{recur_equ}) becomes
\begin{equation}
	\bar{P}_{k|k-1}^{-1} (\hat{\mathcal{X}}_{k|k} - \mathcal{\hat{X}}_{k|k-1}) = H_{k}^{T} \bar{R}_{k}^{-1} (\mathcal{Y}_{k} - H_{k}\hat{\mathcal{X}}_{k|k}). \label{post_st_equ}
\end{equation}
By solving (\ref{post_st_equ}), we will get
\begin{equation}
	\hat{\mathcal{X}}_{k|k} = \mathcal{\hat{X}}_{k|k-1} + K_{k} (\mathcal{Y}_{k} - H_{k}\mathcal{\hat{X}}_{k|k-1}), \label{opt_value}
\end{equation}
where \begin{equation}
	K_{k} = (\bar{P}_{k|k-1} + H_{k}^T \bar{R}_{k}^{-1} H_{k})^{-1} H_{k}^T \bar{R}_{k}^{-1}. \label{K_pre}
\end{equation}

Now, applying Sherman-Morrison-Woodbury matrix identity \cite{riedel1992sherman} in (\ref{K_pre}), (\ref{K_value}) can be obtained.
Using (\ref{Post_state}), posterior error covariance can be calculated as
\begin{equation}
	\begin{split}
		P_{k|k} = &~ \mathbf{E}[(\mathcal{X}_{k} - \mathcal{\hat{X}}_{k|k})(\mathcal{X}_{k} - \mathcal{\hat{X}}_{k|k})^T] \\
		= &~ \mathbf{E} [((\mathcal{X}_{k} - \mathcal{\hat{X}}_{k|k-1}) - K_{k} (H_{k} \mathcal{X}_{k} + r_{k} - H_{k}\mathcal{\hat{X}}_{k|k-1})) ((\mathcal{X}_{k} - \mathcal{\hat{X}}_{k|k-1}) - K_{k} (H_{k} \mathcal{X}_{k} + r_{k} - H_{k}\mathcal{\hat{X}}_{k|k-1}))^T] \\
		= &~ \mathbf{E} [((\mathcal{X}_{k} - \mathcal{\hat{X}}_{k|k-1}) - K_{k} H_{k} (\mathcal{X}_{k} - \mathcal{\hat{X}}_{k|k-1}) - K_{k} r_{k} ) ((\mathcal{X}_{k} - \mathcal{\hat{X}}_{k|k-1}) - K_{k} H_{k} (\mathcal{X}_{k} - \mathcal{\hat{X}}_{k|k-1}) - K_{k} r_{k} )^T] \\
		%= &~ \mathbf{E}[(\mathcal{X}_{k} - \mathcal{\hat{X}}_{k|k-1})(\mathcal{X}_{k} - \mathcal{\hat{X}}_{k|k-1})^T] - \mathbf{E}[ (\mathcal{X}_{k} - \mathcal{\hat{X}}_{k|k-1}) (\mathcal{X}_{k} - \mathcal{\hat{X}}_{k|k-1})^T H_{k}^T K_{k}^T] - \mathbf{E} [K_{k} H_{k} (\mathcal{X}_{k} - \mathcal{\hat{X}}_{k|k-1}) (\mathcal{X}_{k} - \mathcal{\hat{X}}_{k|k-1})^T] \\ 
		%&~ + \mathbf{E} [K_{k} H_{k} (\mathcal{X}_{k} - \mathcal{\hat{X}}_{k|k-1}) (\mathcal{X}_{k} - \mathcal{\hat{X}}_{k|k-1})^T H_{k}^T K_{k}^T] + \mathbf{E} [K_{k} v_{k} v_{k}^T K_{k}^T]\\
		= &~ P_{k|k-1} - P_{k|k-1} H_{k}^T K_{k}^T - K_{k} H_{k} P_{k|k-1} + K_{k} H_{k} P_{k|k-1} H_{k}^T K_{k}^T + K_{k} R_{k} K_{k}^T \\
		= &~ (I - K_{k} H_{k}) P_{k|k-1} (I - K_{k} H_{k})^T + K_{k} R_{k} K_{k}^T.
	\end{split} \label{post_err_cov}
\end{equation} 	
\begin{flushright}
	$\blacksquare$	
\end{flushright}
%\end{pf}

\begin{algorithm} \label{algo_1}
	\caption{Fixed point iteration to calculate $\mathcal{\hat{X}}_{k|k}$}
	\begin{algorithmic}
		%\begin{equation*}
		\STATE $$[\mathcal{\hat{X}}_{k|k}] := FPI[\mathcal{\hat{X}}_{k|k}, \mathcal{\hat{X}}_{k|k-1}, \mathcal{Y}_{k}]$$
		%\end{equation*}
		\begin{enumerate}
			\item Set initial iteration $t = 0$, select $t_{max}$ value and $\hat{\mathcal{X}}_{k|k}^{(t)} = \hat{\mathcal{X}}_{k|k-1}$. \\
			\item \textbf{for} $t=1:t_{max}$
			\begin{itemize}
				\item Calculate $e_{p,k} = - B_{p,k|k-1}^{-1} (\mathcal{\hat{X}}_{k|k}^{(t)} - \mathcal{\hat{X}}_{k|k-1})$ and $e_{r,k} = B_{r,k}^{-1} (\mathcal{Y}_{k} - H_{k}\mathcal{\hat{X}}_{k|k}^{(t)})$.
				\item Calculate $\Pi_{p,k}$ and $\Pi_{r,k}$.
				\item Calculate $K_{k}$ by (\ref{K_value}).
				\item Calculate $\hat{\mathcal{X}}_{k|k}^{(t+1)}$ by (\ref{Post_state}).
				\item \textbf{if} $\frac{||\mathcal{\hat{X}}_{k|k}^{(t+1)} - \mathcal{\hat{X}}_{k|k}^{(t)}||}{||\mathcal{\hat{X}}_{k|k}^{(t)}||} \leq \epsilon$, where $\epsilon$ is the threshold value
				\begin{itemize}
					\item break and update $\mathcal{\hat{X}}_{k|k} = \mathcal{\hat{X}}_{k|k}^{(t+1)}$.
				\end{itemize}
				\textbf{else}  
				\begin{itemize}
					\item $t = t + 1$ and continue iteration.
				\end{itemize} 
				\item \textbf{end if}
			\end{itemize}
			\item \textbf{end for}
		\end{enumerate}
	\end{algorithmic}
\end{algorithm}

\begin{algorithm} \label{algo_2}
	\caption{}
	\begin{algorithmic}
		\STATE {}
		\begin{enumerate}
			\item Set initial values of $\mathcal{\hat{X}}_{0|0}$ and $P_{0|0}$.
			\item Calculate $\mathcal{\hat{X}}_{k|k-1}$ and $P_{k|k-1}$ by using (\ref{Pri_state}) and (\ref{Pri_cov}).
			\item Calculate $B_{p,k|k-1}$ and $B_{r,k}$ using (\ref{v_k_cov}).
			\item Calculate posterior state: \\
			$\mathcal{\hat{X}}_{k|k} = FPI[\mathcal{\hat{X}}_{k|k}, \mathcal{\hat{X}}_{k|k-1}, \mathcal{Y}_{k}]$
			\item Calculate $P_{k|k}$ by (\ref{post_err_cov}).
		\end{enumerate}
	\end{algorithmic}
\end{algorithm}

\begin{remark}
	The error matrix $e_{k}$ plays an important role in our filtering algorithm. It can be observed that the matrix $\Pi_{k} = \begin{bmatrix}
		\Pi_{p,k} & 0 \\
		0 & \Pi_{r,k}
	\end{bmatrix}$ contains the elements of the error matrix which is unavailable to us. Because in practical scenario we don't have the access to the true states, this error matrix can not be calculated directly. Rather we need to adopt some approximate value that can be calculated using the iterative method explained in algorithm 1. It is interesting to note that we are using posterior estimate to calculate the error and to calculate posterior estimate, we need error matrix. Hence this becomes a fixed point iteration and we choose an initial value of posterior estimate as explained in algorithm 1.
\end{remark}
\begin{remark}
	In (\ref{diff_J}), $\rho_{p,i}$ and $\rho_{r,i}$ denote the weighted past errors. It is obvious that when we are calculating the estimation at $k^{th}$ time step, the past errors are already optimized. Also due to the lower value of $\mu_{1}$, the weighted past errors are very less as compare to weighted present error. Hence , it can be ignored as compare to the weighted present error.
\end{remark}	

\begin{remark} \label{rem:3}
	The proposed filter is very sensitive to the kernel bandwidth $\sigma$. For $\sigma \to \infty$, $\Pi_{p,k} = I_{n}$, $\Pi_{r,k} = I_{m}$ and the proposed RMCKF becomes RSKF. For lower value of $\sigma$, the algorithm may not converge. This is explained in upcoming section.
\end{remark}

%\subsection{Calculation of Error Matrix}

%\subsection{Selection of Risk Sensitive Parameter}
\begin{remark}
	The risk sensitive parameter $\mu_{1}$ acts as tuning parameter in this algorithm. Increasing the value of $\mu_{1}$ will increase the robustness of the filtering technique. However during choosing the parameter value, we have to ensure that the condition $(P_{k-1|k-1}^{-1} - 2 \mu_{1} I) > 0$ should be satisfied in order to keep the error covariance matrix positive definite at each time step. Further selection of $\mu_{2}$ should be such that $\Pi_{p,k}$ and $\Pi_{r,k}$ should not be singular at any propagation step. It is interesting to note that although we consider $\mu_1$ and $\mu_2$ as constant, they can also vary with time provided the above condition is satisfied at each time step. 
\end{remark}

\subsection{Selection of Kernel Bandwidth $(\sigma)$}
Kernel bandwidth $(\sigma)$ is an important and sensitive parameter and its proper selection is important for an accurate estimation. There is no reliable method available in literature to get the optimal value of kernel bandwidth. Though a few papers \cite{cinar2012hidden, hou2017maximum, fakoorian2019maximum} addressed the challenge regarding kernel bandwidth selection and proposed some measurement based equation to find the best value at each time step. In \cite{cinar2012hidden}, the norm of the observation error is heuristically considered as the kernel bandwidth value \emph{i.e.} $\sigma_{k} = ||\mathcal{Y}_{k} - H_{k} \mathcal{\hat{X}}_{k}||$ which denotes the Euclidean distance between actual measurement and estimated measurement. In \cite{hou2017maximum}, rather than Euclidean distance, the authors considered Mahalanobis distance that is $\sigma_{k} = ||(\mathcal{Y}_{k} - H_{k} \mathcal{\hat{X}}_{k})^T R_{k}^{-1} (\mathcal{Y}_{k} - H_{k} \mathcal{\hat{X}}_{k})||$. In \cite{fakoorian2019maximum}, the authors consider sum of weighted innovation and weighted error covariance and they took $\sigma_{k} = (||r_{k}||_{R_{k}^{-1}} + H_{k} P_{k|k-1}H_{k}^T)^{-1}$ where $r_{k} = (\mathcal{Y}_{k} - H_{k} \mathcal{\hat{X}}_{k})$. However all the above methods don't guarantee the optimal solution. \\
In this paper, we propose an error based cost function consisting the correntropy criteria inspired by the Eqn.(9) of \cite{singh2010kernel} to find out the kernel bandwidth value at each time step of estimation. We define our cost criteria as 
\begin{equation}
	\begin{split}
		J_{KB}(e_{k,i}) 
		&~ = log(\frac{1}{L} \sum_{i=1}^{L} G_{\sigma}(e_{k,i})) \\
		&~ = log(\frac{1}{L} \sum_{i=1}^{L} \text{exp} (-\frac{e_{k,i}^2}{2 \sigma_{c}^2})). \label{kb_cc} 
	\end{split}
\end{equation}
There is a fundamental dissimilarity in the basic structure of two cost functions defined in (\ref{RSMCC_cost_criteria}) and (\ref{kb_cc}) respectively. 
%In (\ref{RSMCC_cost_criteria}), we have considered process error as well as measurement error which make the dimension of the error matrix ($e_{k}$) $L = n + m$ as explained in section \ref{sec:cf}. Whereas in (\ref{kb_cc}), the cost criteria is defined considering the process error only that leads the dimension of error matrix to be $n$. The (\ref{kb_cc}) can be re-written as
%\begin{equation}
%	J_{KB}(e_{p,k,i}) = log(\frac{1}{n} \sum_{i=1}^{n} \text{exp} (-\frac{e_{p,k,i}^2}{2 \sigma_{c}^2})) \label{kb_cc_exp}.
%\end{equation}
As we have introduced the correntropy function in this cost criteria, there occurs another bandwidth parameter $\sigma_{c}$ in (\ref{kb_cc}). But this $\sigma_{c}$ is different from our earlier kernel bandwidth $\sigma$ and it is a constant value. To find out the $\sigma$ value at $k^{th}$ step, our goal will be to maximize $J_{KB}(e_{k,i})$ \emph{i.e.}
\begin{equation}
	\sigma_{k}^{*} = \text{arg} \text{ $\displaystyle \max_{\sigma_{k}}$} J_{KB}(e_{k,i}) \label{kb_cc_arg}.
\end{equation}
It can be seen that (\ref{kb_cc_arg}) actually symbolises the minimization of error. Hence the $\sigma_{c}$ should be a constant value for the proper identification of the optimal value of $\sigma_{k}$. \\
To identify the desired kernel bandwidth, a numerical search rule is considered. At first we define a range of possible $\sigma$ values and calculate the error at each time step as explained in algorithm 1. Using the error value, $J_{KB}(e_{k,i})$ is obtained for each $\sigma_k$. Further we compare the $J_{KB}(e_{k,i})$ values and pick the maximum one and its corresponding $\sigma_k$. We repeat the same process for every time step.

\section{Convergence and Stability analysis}
In this section, we prove the convergence and stability of the proposed algorithm. Further, convergence of fixed point iteration algorithm mentioned in algorithm-1 is also established. 
\subsection{Convergence of filter}
\begin{lemma}
	If $\tilde{\theta}_{\bar{n}} = \mathcal{X}_{\bar{n}} - \mathcal{\hat{X}}_{\bar{n}|\bar{n}}$ defines the error of the state, then
	$\mathbf{E}[||\tilde{\theta}_{\bar{n}}||^2] \leq \mathbf{E}[||S_{\bar{n}|\bar{n}}||^2] \mathbf{E}[V_{\bar{n}}]$, where $V_{\bar{n}}$ is a Lyapunov function defined by $V_{\bar{n}} = \tilde{\theta}_{\bar{n}}^T P_{\bar{n}|\bar{n}}^{-1} \tilde{\theta}_{\bar{n}}$, $P_{\bar{n}|\bar{n}} = S_{\bar{n}|\bar{n}} S_{\bar{n}|\bar{n}}^T$ and $\bar{n}$ is the maximum time step satisfying the condition $\bar{n} \geq k \geq 0$.
\end{lemma}

\textbf{Proof:}
Recalling the matrix reversal law, we can write $S_{\bar{n}|\bar{n}} S_{\bar{n}|\bar{n}}^{-1} = 1$.
Hence, 
\begin{equation}
	\begin{split}
		\mathbf{E}[||\tilde{\theta}_{\bar{n}}||^2] 
		&~ = \mathbf{E}[||S_{\bar{n}|\bar{n}} S_{\bar{n}|\bar{n}}^{-1} \tilde{\theta}_{\bar{n}}||^2] \\
		&~ \leq \mathbf{E}[||S_{\bar{n}|\bar{n}}||^2 ||S_{\bar{n}|\bar{n}}^{-1} \tilde{\theta}_{\bar{n}}||^2] \\
		&~ \leq \mathbf{E}[||S_{\bar{n}|\bar{n}}||^2] \mathbf{E} [||S_{\bar{n}|\bar{n}}^{-1} \tilde{\theta}_{\bar{n}}||^2] \\
	\end{split}
\end{equation}
Using the matrix norm property $||A||^2 = ||A^T A|| = || A A^T||$ $\forall A$, we get
\begin{equation}
	\begin{split}
		\mathbf{E}[||\tilde{\theta}_{\bar{n}}||^2] 
		&~ \leq \mathbf{E}[||S_{\bar{n}|\bar{n}}||^2] \mathbf{E} [||\tilde{\theta}_{\bar{n}}^T S_{\bar{n}|\bar{n}}^{-1T} S_{\bar{n}|\bar{n}}^{-1} \tilde{\theta}_{\bar{n}}||] \\
		&~ \leq \mathbf{E}[||S_{\bar{n}|\bar{n}}||^2] \mathbf{E} [||\tilde{\theta}_{\bar{n}}^T P_{\bar{n}|\bar{n}}^{-1} \tilde{\theta}_{\bar{n}}||] \\
		&~ \leq \mathbf{E}[||S_{\bar{n}|\bar{n}}||^2] \mathbf{E}[V_{\bar{n}}]
	\end{split}
\end{equation}

\begin{lemma}
	The filter is convergent if  for $\bar{n} \to \infty$, the following condition satisfies
	\begin{equation*}
		\mathbf{E}[||\tilde{\theta}_{\bar{n}}||] \to 0
	\end{equation*}
\end{lemma}
\textbf{Proof:} 
(See Appendix A).

\begin{remark}
	The notation $O$ defines the Landau’s Symbol (also called big O notation). It indicates the rate of how fast or slow a function will decay or grow. More details can be found in Appendix B of \cite{9452347}.
\end{remark}

\subsection{Convergence of fixed point iteration}

\begin{lemma}
	The convergence of the fixed point iteration is guaranteed for $\beta > 0$ and $0 < \alpha < 1$, where  the initial vector $||\hat{\mathcal{X}}_{0|0}||_{1} \leq \beta$ and $\forall$ $\hat{\mathcal{X}}_{k|k}$, the following conditions hold
	
	% $\mathcal{X}_{k} \in \{\mathcal{X}_{k} \in \mathbb{R}^{n}: ||\mathcal{X}_{k}||_{1} \leq \beta\}$
	
	\begin{align} 
		&(i) ~~~~ 
		||\hat{f}(\hat{\mathcal{X}}_{k|k})||_{1} \leq \beta \label{cond_1} , \\
		& (ii) ~~~~
		||\frac{\partial \hat{f}(\hat{\mathcal{X}}_{k|k})}{\partial \mathcal{X}_{k}}||_{1} \leq \alpha \label{cond_2} ,
	\end{align}
	where $\hat{f}(\hat{\mathcal{X}}_{k|k})$ denotes the fixed point iteration such that $\hat{f}(\hat{\mathcal{X}}_{k|k}) = (W_{k}^T \Pi_{k} W_{k})^{-1} W_{k} \Pi_{k} D_{k}$.
\end{lemma}

%\begin{pf}
\textbf{Proof:}
(See Appendix B).
%\end{pf}

\subsection{Stability Analysis}
We assume the system parameters $F_{k}$, $H_{k}$, $Q_{k}$ and $R_{k}$ are stochastically bounded and the system uncertainty parameter $\Delta F_{k}$ as defined in (\ref{process_equ}) is finite. We consider that the equivalent covariance of measurement noise $R_{k}$ is non-zero and finite so that $R_{k}^{-1}$ is non-singular and bounded. Let us define the controllability Grammian matrix and observability Grammian matrix for the system defined in (\ref{process_equ}), (\ref{measure_equ}) \cite{jazwinski2007stochastic} as 
\begin{equation}
	\mathcal{C}_{k,k-l} = \sum_{i=k-l}^{k-1} (F_{k,i+1} + \Delta F_{k,i+1})^T Q_{i} (F_{k,i+1} + \Delta F_{k,i+1}) \label{con_mat},
\end{equation}
\begin{equation}
	\mathcal{O}_{k,k-l} = \sum_{i=k-l}^{k} (F_{i,k} + \Delta F_{i,k})^T H_{i}^T R_{i}^{-1} H_{i} (F_{i,k} + \Delta F_{i,k}) \label{obs_mat},
\end{equation}
where $l$ is a positive integer. Now, the system (\ref{process_equ}), (\ref{measure_equ}) is uniformly completely observable and uniformly completely controllable if the observability Grammian matrix $\mathcal{O}_{k,k-l}$ and controllability Grammian matrix $\mathcal{C}_{k,k-l}$ are finite and bounded, \emph{i.e.} $0< k_{1} I < \mathcal{O}_{k,k-l} < k_{2} I$ and $0< k_{3} I < \mathcal{C}_{k,k-l} < k_{4} I$ where $k_{1}$, $k_{2}$,
$k_{3}$ and $k_{4}$ are real and positive. We consider that the equivalent posterior error covariance $P_{k|k}$ is positive definite for any $k$. 

\begin{lemma}
	The filter is stable if the system (\ref{process_equ}), (\ref{measure_equ}) is uniformly completely observable and uniformly completely controllable and if $P_{0} > 0$ and $(P_{k-1|k-1}^{-1} - 2 \mu_{1}I) > 0$ for all $k>l$ provided that $- I < \mathcal{O}_{k,k-l}^{-1} \Delta \mathcal{O}_{k,k-l} < I$, where
	$\Delta \mathcal{O}_{k,k-l} = \sum_{i = k-l}^{k} F_{i,k}^T H_{i}^T R_{i}^{-1} H_{i} \Delta F_{i,k}$.
\end{lemma}

\textbf{Proof:} 
Considering the assumptions stated above and following the Appendix B of \cite{tiwari2022risk}, the stability of the filter can be proved.

\begin{flushright}
	$\blacksquare$	
\end{flushright}

\section{Simulation Results}
In this section, we will simulate some numerical examples to verify the effectiveness of our proposed algorithm. We will compare our results with already developed filters available in literature. We will also make a comparative study of our results for different values of the uncertainty parameter. 
\subsection{Problem 1}
Let us consider a system modelled as 
\begin{equation*}
	\mathcal{X}_{k+1} = (F + \Delta F) \mathcal{X}_{k} + G q_k, 
\end{equation*}
\begin{equation*}
	\mathcal{Y}_{k} = H \mathcal{X}_{k} + r_k,
\end{equation*}
where $F = \begin{bmatrix}
	0.99 & 0.01 \\
	0 & 0.99
\end{bmatrix}$, $\Delta F = \begin{bmatrix}
	0 & \delta \\
	0 & 0
\end{bmatrix}$, $G = \begin{bmatrix}
	5 \\ 1
\end{bmatrix}$ and $H = \begin{bmatrix}
	1 & -1
\end{bmatrix}$. $\delta$ is the uncertainty in the given system which is bounded and does not effect the system's stability. We consider $| \delta| \leq 0.5$. The process noise $q_{k}$ and measurement noise $r_{k}$ are zero-mean and non-Gaussian in nature which are modelled as sum of Gaussian distributions. We consider $q_{k} = 0.8 \mathcal{N}(0,0.01) + 0.2 \mathcal{N}(0, 1)$ and $r_{k} = 0.8 \mathcal{N}(0, 1) + 0.2\mathcal{N}(0,1000)$. The initial state of truth is taken as $x_{0} = \begin{bmatrix}
	10 & 20
\end{bmatrix}^T$ and initial error covariance is $P_{0} = \begin{bmatrix}
	35^2 & 0 \\
	0 & 70^2
\end{bmatrix}$. For estimation we have considered random initial state with mean $x_{0}$ and covariance $P_{0}$.
The risk parameters $\mu_{1,k-1}$ is selected in such a way that $(P_{k-1|k-1}^{-1} - 2 \mu_{1,k-1}I) > 0$ satisfies at each time step $k$ and $\mu_{2, k-1}$ is an arbitrary scalar value. We have considered equivalent covariance of process noise and measurement noise for filtering.

\begin{figure}
	\begin{center}
		\includegraphics[width=12cm]{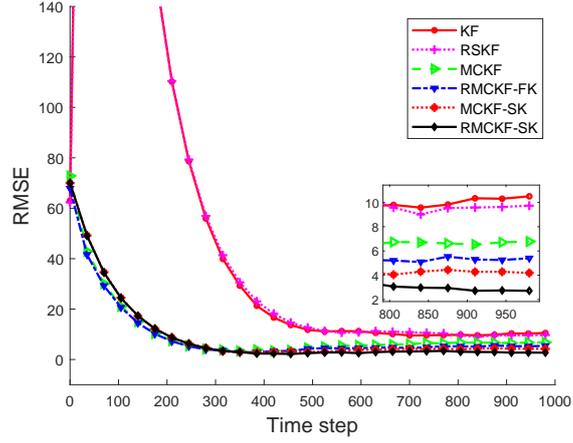}
		\caption{RMSE plot} 
		\label{fig:prb_1_rmse_plot}
	\end{center}
\end{figure}

\begin{figure}
		\begin{center}
		\includegraphics[width=12cm]{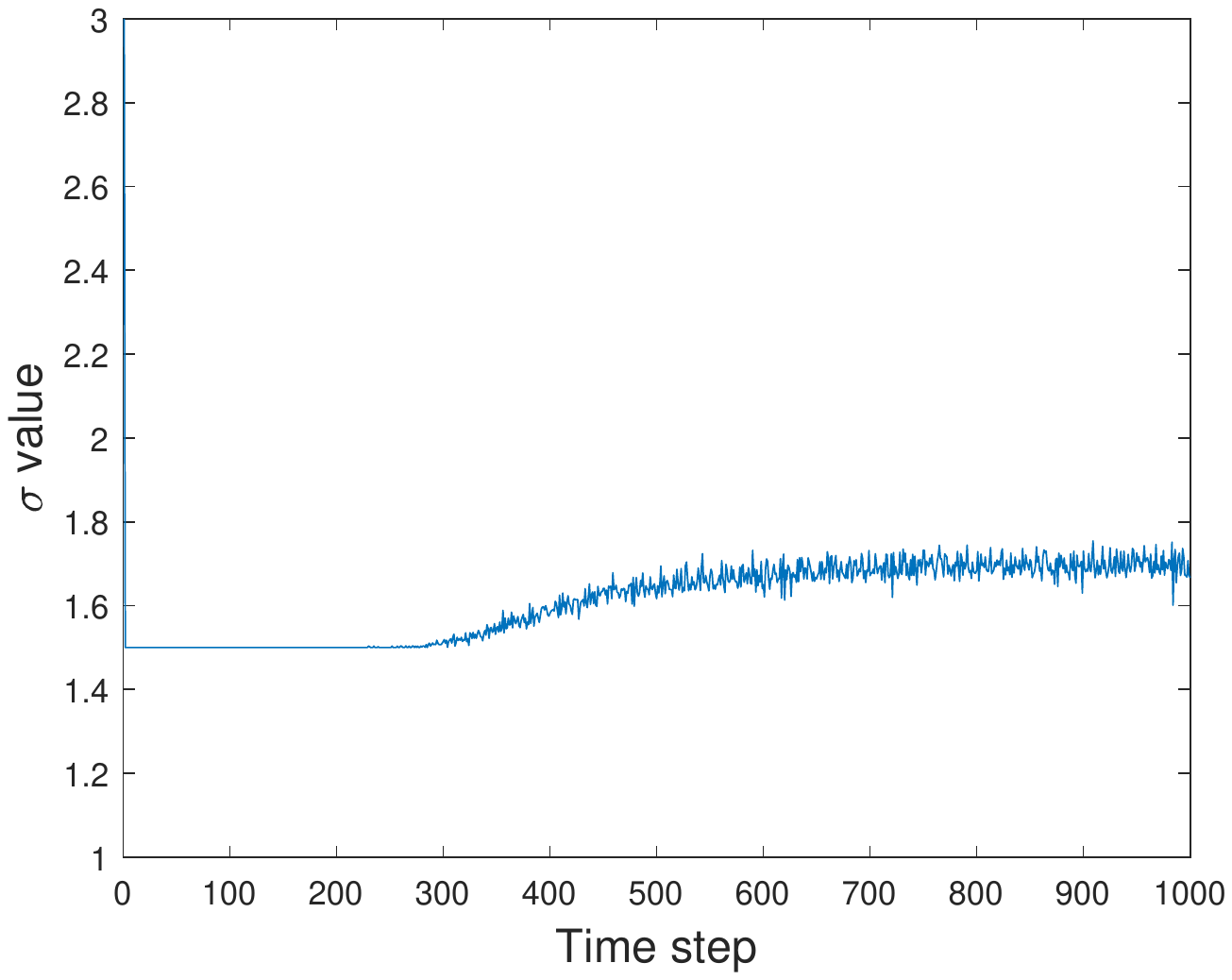}
		\caption{Kernel Bandwidth $(\sigma)$ value at each time step} 
		\label{fig:prb_1_kb_selection}
		\end{center}
\end{figure}

Fig.\ref{fig:prb_1_rmse_plot} compares the root mean square error (RMSE) of state 2 for Kalman filter (KF), risk sensitive Kalman filter (RSKF), maximum correntropy Kalman filter (MCKF), robust maximum correntropy Kalman filter with fixed kernel bandwidth (RMCKF-FK), maximum correntropy Kalman filter with selected kernel bandwidth  (MCKF-SK), and robust maximum correntropy Kalman filter with selected kernel bandwidth (RMCKF-SK). The same for state 1 is not shown due to similar characteristics. It can be seen that RMCKF-SK gives better result as compare to the other filters. In Fig.\ref{fig:prb_1_kb_selection} the selected values of kernel bandwidth$(\sigma)$ for each time step $k$ is shown. 

\begin{table}[htbp] \label{tab:prb_1_rmse}
	\centering
	\caption{Comparision of avg-RMSE for different filters}  
	\begin{tabular}{ |p{0.4cm}|p{0.6cm}|p{0.7cm}|p{0.7cm}|p{1.1cm}|p{0.9cm}|p{1.1cm}| }
		\cline{1-7}
		$\delta$ & KF  & RSKF & MCKF & RMCKF-FK & MCKF-SK & RMCKF-SK            \\
		\cline{1-7}
		0	&	0.55	&	0.59	&	0.52	&	0.55	&	0.43	&	0.50	\\
		\cline{1-7}
		0.1	&	2.29	&	2.23	&	2.12	&	2.08	&	2.07	&	1.97	\\
		\cline{1-7}
		0.2	&	4.45	&	4.26	&	3.97	&	4.11	&	3.85	&	3.15	\\
		\cline{1-7}
		0.3	&	6.64	&	6.32	&	5.47	&	5.05	&	4.22	&	3.41	\\
		\cline{1-7}
		0.4	&	8.85	&	8.49	&	5.64	&	4.96	&	4.68	&	2.77	\\
		\cline{1-7}
		0.5	&	10.65	&	10.43	&	6.23	&	4.96	&	4.06	&	2.67	\\
		\cline{1-7}
	\end{tabular}
\end{table}	

A detailed study on the performance of KF, RSKF, MCKF, RMCKF-FK, MCKF-SK and RMCKF-SK is performed and the results are shown in Table-1. It is observed that RMCKF-SK outperforms all other filters. Some observations that can be figured out from the table-1 are as follows:

\begin{enumerate}
	\item When there is no uncertainty in the system model \emph{i.e.} $\delta = 0$, robust filters are less accurate as compared to normal filters. It suggests that at $\delta = 0$, KF, MCKF and MCKF-SK are better than RSKF, RMCKF and RMCKF-SK respectively. 
	\item With the increase of uncertainty parameter $\delta$, rmse also increases for all the filters. moreover it is notable that RMCKF-FK and RMCKF-SK are always better than MCKF and MCKF-SK for non-zero $\delta$. 
	\item MCKF-SK and RMCKF-SK are always better than MCKF and RMCKF-FK respectively. It indicates that kernel bandwidth is a sensitive and very important parameter in correntropy based filters. Proper selection of kernel bandwidth always provides better results.
\end{enumerate}

\subsection{Problem 2}
Let us consider a system moving with constant acceleration. Also we consider that due to some external force or some internal disturbance, the model have uncertainty. Define the states as $\mathcal{X}_{k} = \begin{bmatrix}
	s_k & v_k & a_k
\end{bmatrix}^T$, where $s_k$, $v_k$ and $a_k$ denotes position, velocity and acceleration of the system respectively. In discrete time, the system can be modelled as
\begin{equation*}
	\mathcal{X}_{k+1} = (F + \Delta F) \mathcal{X}_{k} + q_k, 
\end{equation*}
\begin{equation*}
	\mathcal{Y}_{k} = H \mathcal{X}_{k} + r_k,
\end{equation*}
where 
\begin{equation*}
	F = \begin{bmatrix}
		1 & T & \frac{1}{2}T^2 \\
		0 & 1 & T \\
		0 & 0 & 1
	\end{bmatrix},
\end{equation*}
\begin{equation*} \Delta F =
	\begin{bmatrix}
		0 & 0 & \delta_{1} T^2 \\
		0 & 0 & \delta_{2} T \\
		0 & 0 & 0
	\end{bmatrix}
\end{equation*}

$H = \begin{bmatrix}
	1 & 1 & 0
\end{bmatrix}$, the sampling time T is considered 0.1 min, and $q_{k} = \begin{bmatrix}
	q_{1k} & q_{2k} & q_{3k}
\end{bmatrix}^T$ and $|\delta_{1}| \leq 0.005$ $\&$ $|\delta_{2}| \leq 0.05$. The term $\delta_{1} T^2 a_{k}$ represents an uncertainty in position that is modelled as the function of acceleration $(a_{k})$ and sampling time $(T)$. Also the term $\delta_{2} T a_{k}$ defines the uncertainty in velocity as the function of acceleration $(a_{k})$ and sampling time $(T)$. We have considered the uncorrelated noises in Gaussian mixture form with zero-mean, distributed as $q_{1k} = q_{2k} = q_{3k} = 0.9\mathcal{N}(0,0.0005) + 0.1\mathcal{N}(0,0.05)$ and $r_k = 0.8\mathcal{N}(0,0.005) + 0.2\mathcal{N}(0,50)$. The initial state $x_0 = \begin{bmatrix}
	50 & 4 & 1
\end{bmatrix}^T$ and initial error covariance matrix is taken as $P_{0} = \begin{bmatrix}
	0.5^2 & 0 & 0 \\
	0 & 0.5^2 & 0 \\
	0 & 0 & 0.1^2
\end{bmatrix}$.  
The risk sensitive parameters $\mu_{1,k-1}$ and $\mu_{2,k-1}$ are selected as explained problem 1.

\begin{figure}
	\begin{center}
		\includegraphics[width=12cm]{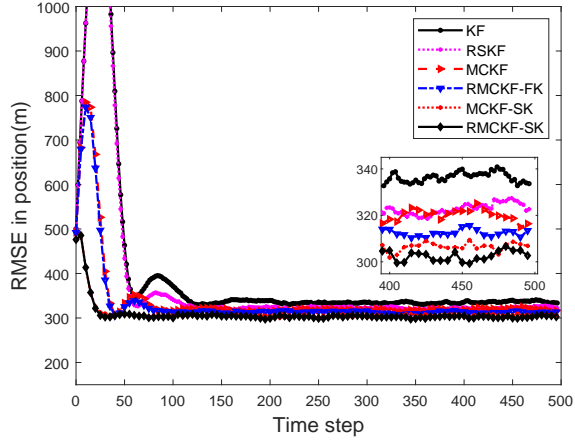}
		\caption{RMSE in position} 
		\label{fig:prob2_rmse_pos}
	\end{center}
\end{figure}

\begin{figure}
	\begin{center}
		\includegraphics[width=12cm]{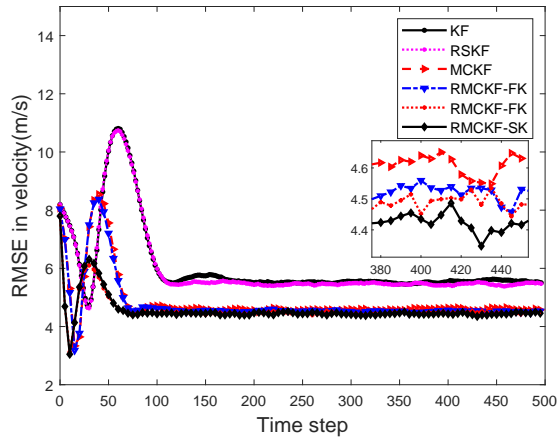}
		\caption{RMSE in velocity} 
		\label{fig:prob2_rmse_vel}
	\end{center}
\end{figure}

\begin{figure}
	\begin{center}
		\includegraphics[width=12cm]{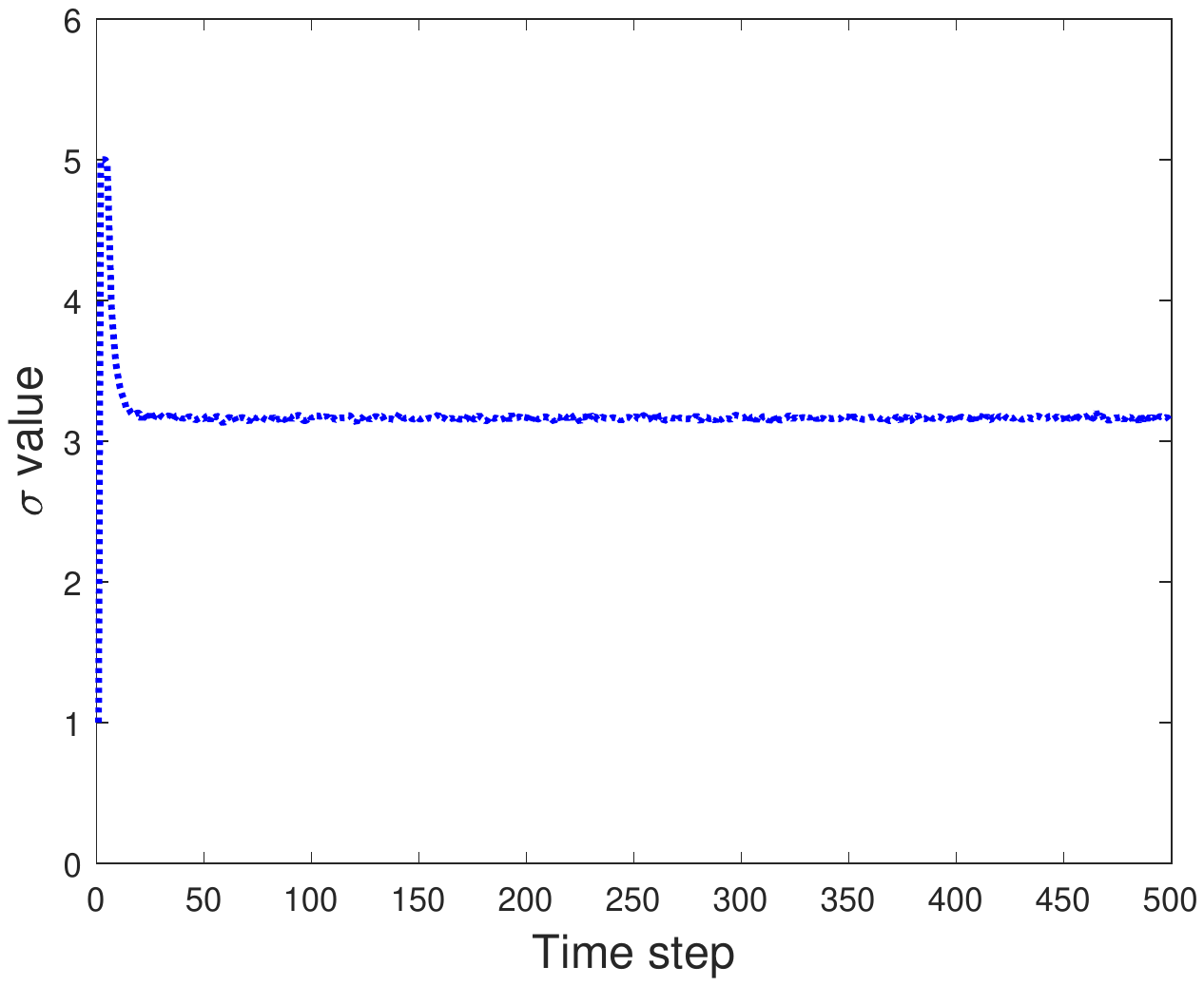}
		\caption{Kernel Bandwidth $(\sigma)$ value at each time step} 
		\label{fig:prob2_sigma}
	\end{center}
\end{figure}

\begin{table*}[t] \label{tab:prb_2_rmse_pos}
	\centering
	\caption{Comparision of avg-RMSE in position and velocity (m/s) (m) for different filters}  
	\begin{tabular}{ |p{0.6cm}|p{0.8cm}|p{0.8cm}|p{0.8cm}|p{1.2cm}|p{1.0cm}|p{1.2cm}|p{0.8cm}|p{0.8cm}|p{0.8cm}|p{1.2cm}|p{1.0cm}|p{1.2cm}| }
		\cline{1-13}
		{} & \multicolumn{6}{c|}{Position} & \multicolumn{6}{|c|}{Velocity} \\
		\cline{2-13}
		$\delta_2$ & KF  & RSKF & MCKF & RMCKF-FK & MCKF-SK & RMCKF-SK   & KF  & RSKF & MCKF & RMCKF-FK & MCKF-SK & RMCKF-SK     \\
		\cline{1-13}
		0.000	&	407.10	&	409.86	&	341.67	&	343.23	&	305.76	&	312.56 &	5.88	&	5.91	&	4.75	&	4.79	&	4.42	&	4.57	\\
		\cline{1-13}
		0.005	&	406.47	&	409.50	&	338.90	&	339.63	&	306.56	&	315.00	&	5.88	&	5.91	&	4.76	&	4.75	&	4.47	&	4.61 \\
		\cline{1-13}
		0.01	&	409.08	&	407.68	&	340.44	&	336.87	&	307.38	&	306.67	&	5.99	&	5.93	&	4.79	&	4.70	&	4.44	&	4.33 \\
		\cline{1-13}
		0.015	&	405.33	&	405.09	&	341.10	&	337.07	&	307.33	&	300.37	&	5.95	&	5.91	&	4.79	&	4.75	&	4.45	&	4.36 \\
		\cline{1-13}
		0.02	&	405.01	&	403.26	&	343.62	&	333.84	&	307.89	&	304.54	&	5.99	&	5.90	&	4.80	&	4.73	&	4.50	&	4.38 \\
		\cline{1-13}
		0.025	&	405.67	&	401.03	&	342.56	&	333.11	&	306.56	&	302.28	&	5.99	&	5.90	&	4.80	&	4.60	&	4.48	&	4.43 \\
		\cline{1-13}
		0.03	&	407.21	&	402.73	&	341.89	&	333.57	&	308.13	&	303.05	&	6.00	&	5.95	&	4.85	&	4.75	&	4.52	&	4.42 \\
		\cline{1-13}
		0.035	&	405.60	&	403.01	&	342.49	&	327.74	&	307.90	&	304.79	&	6.01	&	5.94	&	4.82	&	4.73	&	4.52	&	4.47 \\
		\cline{1-13}
		0.04	&	404.34	&	401.90	&	342.55	&	337.27	&	307.81	&	300.17	&	6.01	&	5.89	&	4.85	&	4.80	&	4.52	&	4.43 \\
		\cline{1-13}
		0.045	&	404.75	&	401.61	&	343.53	&	332.85	&	307.64	&	304.70	&	6.07	&	5.96	&	4.89	&	4.74	&	4.55	&	4.50 \\
		\cline{1-13}
		0.05	&	399.94	&	397.95	&	343.57	&	330.89	&	308.66	&	304.95	&	6.04	&	6.00	&	4.87	&	4.74	&	4.53	&	4.44 \\
		\cline{1-13}
	\end{tabular}
\end{table*}

The root mean square error (RMSE) of position and velocity are shown respectively in Fig.\ref{fig:prob2_rmse_pos} and Fig.\ref{fig:prob2_rmse_vel} for KF, RSKF, MCKF, RMCKF-FK, MCKF-SK and RMCKF-SK. It is observed that RMCKF-SK provides better result as compared to other filters. The selected values of kernel bandwidth is shown in Fig.\ref{fig:prob2_sigma}. Table.2 shows the variation in position RMSE and velocity RMSE respectively with the change in system uncertainty parameter $\delta_{2}$. It can be concluded from the table that RMCKF-SK gives better result as compare to all other mentioned filters. It may be questionable that why do we vary $\delta_{2}$? From the system description, it can arguably said that $\delta_{1}$ impact the position only, whereas $\delta_{2}$ directly impact velocity and hence position is also getting impacted. So, varying $\delta_{2}$ means the changing the uncertainty in both position and velocity. That's why we choose to vary $\delta_{2}$.

\section{Conclusion}
We have developed a new filtering algorithm to deal with uncertain system model in presence of non-Gaussian noises where nominal robust Kalman filter fails. We have proposed a new cost function using maximum correntropy criteria and by maximizing this, our proposed filtering recursion equations are established. We also presented a new numerical approach to select kernel bandwidth at each time step for better performance. The condition of stability, convergence of filter and convergence of fixed point iteration to calculate posterior state is presented

\bibliographystyle{ieeetr}
\bibliography{ref_paper1.bib}

\section*{Appendices}

\subsection*{A. ~ Proof of lemma 3}
From Lemma 2, it can be said that $\mathbf{E}[||\tilde{\theta}_{\bar{n}}||] \to 0$ if $\mathbf{E} [V_{\bar{n}}] \to 0$. So our task reduces to prove that $\mathbf{E} [V_{\bar{n}}] \to 0$ for $\bar{n} \to \infty$.
Considering $(I - K_{k} H_{k}) = G_{k}$ in (\ref{post_err_cov}), and using (\ref{Pri_cov}), we get
\begin{equation}
	P_{k|k} = G_{k} F_{k-1} (P_{k-1|k-1}^{-1} - 2 \mu_{1} I)^{-1} F_{k-1}^T G_{k}^T + K_{k} R_{k} K_{k}^T + \tilde{Q}_{k-1}, \label{post_err_cov_recon}
\end{equation}
where $\tilde{Q}_{k-1} = G_{k} Q_{k-1} G_{k}^T$. The error $\tilde{\theta}_{k}$ can be constructed as.

\begin{equation}
	\begin{split}
		\tilde{\theta}_{k} 
		&~ = (\mathcal{X}_{k} - \mathcal{\hat{X}}_{k|k}) \\
		&~ = \mathcal{X}_{k} - \mathcal{\hat{X}}_{k|k-1} - K_{k} (\mathcal{Y}_{k} - H_{k}\mathcal{\hat{X}}_{k|k-1}) \\
		&~ = (\mathcal{X}_{k} - \mathcal{\hat{X}}_{k|k-1}) - K_{k}(H_{k} \mathcal{X}_{k} + r_{k} - H_{k}\mathcal{\hat{X}}_{k|k-1}),
	\end{split}
\end{equation}
Substituting $\mathcal{X}_{k}$ from (\ref{process_equ}) and $\mathcal{\hat{X}}_{k|k-1}$ from (\ref{Pri_state}) and rearranging the terms, we get

\begin{equation}
	\begin{split}
		\tilde{\theta}_{k}
		&~ = G_{k} F_{k-1} \tilde{\theta}_{k-1} + G_{k} q_{k-1} -K_{k} r_{k} \\
		&~ = G_{k} F_{k-1} \tilde{\theta}_{k-1} + z_{k},
	\end{split}
\end{equation}
where $z_{k} = -K_{k} r_{k} + G_{k} q_{k-1}$.
Using the defined Lyapunov function, it can be written as 
\begin{equation}
	\begin{split}
		V_{k} = &~ [G_{k} F_{k-1} \tilde{\theta}_{k-1} + z_{k}]^T P_{k|k}^{-1} [G_{k} F_{k-1} \tilde{\theta}_{k-1} + z_{k}] \\
		= &~ \tilde{\theta}_{k-1}^T F_{k-1}^T G_{k}^T P_{k|k}^{-1} G_{k} F_{k-1} \tilde{\theta}_{k-1} + \tilde{\theta}_{k-1}^T F_{k-1}^T G_{k}^T P_{k|k}^{-1} z_{k} + z_{k}^T P_{k|k}^{-1} G_{k} F_{k-1} \tilde{\theta}_{k-1} + z_{k}^T P_{k|k}^{-1} z_{k}.
	\end{split} \label{lyap_cons}
\end{equation}

Now,
\begin{equation}
	\begin{split}
		F_{k-1}^T G_{k}^T P_{k|k}^{-1} G_{k} F_{k-1} 
		= &~ F_{k-1}^T G_{k}^T [G_{k} F_{k-1} (P_{k-1|k-1}^{-1} - 2 \mu_{1} I)^{-1} F_{k-1}^T G_{k}^T + K_{k} R_{k} K_{k}^T + \tilde{Q}_{k-1}]^{-1} G_{k} F_{k-1} \\
		\leq &~ F_{k-1}^T G_{k}^T [G_{k} F_{k-1} P_{k-1|k-1}   F_{k-1}^T G_{k}^T + K_{k} R_{k} K_{k}^T + \tilde{Q}_{k-1}]^{-1} G_{k} F_{k-1} \\
		\leq &~ [P_{k-1|k-1} + F_{k-1}^{-1} G_{k}^{-1} (K_{k} R_{k} K_{k} + \tilde{Q}_{k-1}) (F_{k-1}^{-1} G_{k}^{-1})^T]^{-1}
	\end{split} \label{lyap_fir_1}
\end{equation}

Using matrix inversion formula presented in eqn. (16) of \cite{henderson1981deriving} in (\ref{lyap_fir_1}), we get
\begin{equation}
	\begin{split}    
		F_{k-1}^T G_{k}^T P_{k|k}^{-1} G_{k} F_{k-1}  
		\leq &~ P_{k-1|k-1}^{-1} - P_{k-1|k-1}^{-1} F_{k-1}^{-1} G_{k}^{-1} [(K_{k} R_{k} K_{k}^T + \tilde{Q}_{k-1})^{-1} + (G_{k}^{-1})^T (F_{k-1}^{-1})^T P_{k-1|k-1}^{-1} F_{k-1}^{-1} G_{k}^{-1}]^{-1} (G_{k}^{-1})^T (F_{k-1}^{-1})^T P_{k-1|k-1}^{-1} \\
		\leq &~ P_{k-1|k-1}^{-1} - [P_{k-1|k-1} + P_{k-1|k-1} F_{k-1}^T G_{k}^T (K_{k} R_{k} K_{k}^T + \tilde{Q}_{k-1})^{-1} G_{k} F_{k-1} P_{k-1|k-1}]^{-1} \\
		\leq &~ S_{k-1|k-1}^{-1 T} (I - [I + S_{k-1|k-1}^T F_{k-1}^T G_{k}^T (K_{k} R_{k} K_{k}^T + \tilde{Q}_{k-1})^{-1} G_{k} F_{k-1} S_{k-1|k-1}]^{-1}) S_{k-1|k-1}^{-1} 
	\end{split} \label{lyap_fir_11}
\end{equation}
Now, using the property $C A B \leq C ||A|| B$  $\forall$ $A>0, B>0, C>0$, we get
\begin{equation}
	\begin{split}   
		F_{k-1}^T G_{k}^T P_{k|k}^{-1} G_{k} F_{k-1}
		\leq &~ S_{k-1|k-1}^{-1 T} (1 - [1 + \frac{||S_{k-1|k-1}^T F_{k-1}^T G_{k}^T G_{k} F_{k-1} S_{k-1|k-1}||}{||(K_{k} R_{k} K_{k}^T + \tilde{Q}_{k-1})||}]^{-1}) S_{k-1|k-1}^{-1} \
	\end{split} \label{lyap_fir_12}
\end{equation}
Using the matrix norm property $||A^T A|| = ||A A^T||$ $\forall$ $A$, we can write
\begin{equation}
	\begin{split} 
		F_{k-1}^T G_{k}^T P_{k|k}^{-1} G_{k} F_{k-1}
		\leq &~ [1 - (1 + ||(K_{k} R_{k} K_{k}^T + \tilde{Q}_{k-1})^{-1} G_{K} F_{k-1} P_{k-1|k-1} F_{k-1}^T G_{k}^T ||)^{-1}] P_{k-1|k-1}^{-1}
	\end{split} \label{lyap_fir_2}
\end{equation}

Using (\ref{post_err_cov_recon}) in (\ref{lyap_fir_2}), we can write
\begin{equation}
	F_{k-1}^T G_{k}^T P_{k|k}^{-1} G_{k} F_{k-1} 
	\leq [1 - (1 + ||(K_{k} R_{k} K_{k}^T + \tilde{Q}_{k-1})^{-1} P_{k|k} ||)^{-1}] P_{k-1|k-1}^{-1} \label{lyap_fir_3}
\end{equation}
Now, $P_{k-1|k-1} + \tilde{Q}_{k-1} \geq P_{k|k}$. Hence. (\ref{lyap_fir_3}) will be
\begin{equation}
	\begin{split}
		F_{k-1}^T G_{k}^T P_{k|k}^{-1} G_{k} F_{k-1} 
		\leq &~ [1 - (1 + ||\tilde{Q}_{k-1}^{-1} (P_{k-1|k-1} + \tilde{Q}_{k-1})||)^{-1}] P_{k-1|k-1}^{-1} \\
		\leq &~ [1 - (2 + ||\tilde{Q}_{k-1}^{-1}|| ||P_{k-1|k-1}||)^{-1}] P_{k-1|k-1}^{-1} \\
		= &~ P_{k-1|k-1}^{-1} - \frac{P_{k-1|k-1}^{-1}}{2 + ||\tilde{Q}_{k-1}^{-1}|| ||P_{k-1|k-1}||}
	\end{split} \label{lyap_fir}
\end{equation}
Using (\ref{lyap_fir}) in (\ref{lyap_cons}), we get
\begin{equation}
	\begin{split}
		V_{k} \leq &~ V_{k-1} - \frac{V_{k-1}}{2 + ||\tilde{Q}_{k-1}^{-1}|| ||P_{k-1|k-1}||} + \tilde{\theta}_{k-1}^T F_{k-1}^T G_{k}^T P_{k|k}^{-1} z_{k} + z_{k}^T P_{k|k}^{-1} G_{k} F_{k-1} \tilde{\theta}_{k-1} + z_{k}^T P_{k|k}^{-1} z_{k} \\
		\leq &~ V_{k-1} - \frac{V_{k-1}}{2 + ||\tilde{Q}_{k-1}^{-1}|| ||P_{k-1|k-1}||} + 2 ||z_{k}^T P_{k|k}^{-1} G_{k} F_{k-1} \tilde{\theta}_{k-1}|| + ||z_{k}^T P_{k|k}^{-1} z_{k}||
	\end{split} \label{lyap_sim}
\end{equation}
Using the elementary inequality $2|ab| \leq a^2 + b^2$, we can write
\begin{equation}
	\begin{split}
		2 ||z_{k}^T P_{k|k}^{-1} G_{k} F_{k-1} \tilde{\theta}_{k-1}|| 
		\leq &~ 2 ||z_{k}^T S_{k|k}^{-1}|| ||S_{k|k}^{-1} G_{k} F_{k-1} \tilde{\theta}_{k-1}|| \\
		\leq &~ ||z_{k}^T P_{k|k}^{-1} z_{k}|| + ||\tilde{\theta}_{k-1}^T F_{k-1}^T G_{k}^T P_{k|k}^{-1}G_{k} F_{k-1} \tilde{\theta}_{k-1}||
	\end{split}
\end{equation}
Using (\ref{lyap_fir}), we get $\tilde{\theta}_{k-1}^T F_{k-1}^T G_{k}^T P_{k|k}^{-1}G_{k} F_{k-1} \tilde{\theta}_{k-1} \leq V_{k-1}$ and substituting $z_{k} = -K_{k} r_{k} + G_{k} q_{k-1}$ we get $||z_{k}^T P_{k|k}^{-1} z_{k}|| \leq || S_{k|k}^{-1} (-K_{k} r_{k} + G_{k} q_{k-1})||^{2} \leq O(||P_{k|k}^{-1}|| (||r_{k}||^{2} + ||q_{k-1}||^{2}))$.
So, from (\ref{lyap_sim}) we get
\begin{equation}
	V_{k} \leq 2[V_{k-1} - \frac{V_{k-1}}{4 + 2 ||\tilde{Q}_{k-1}^{-1}||||P_{k-1|k-1}||}] + O(||P_{k|k}^{-1}|| (||r_{k}||^{2} + ||q_{k-1}||^{2})).
\end{equation}
Now, consider a function $\phi(\bar{n},k)$ such that 
\begin{equation}
	\phi(\bar{n},k) = (1 - \frac{1}{4 + 2 ||\tilde{Q}_{\bar{n}-1}^{-1}||||P_{\bar{n}-1|\bar{n}-1}||}) \phi(\bar{n}-1,k),
\end{equation}
$\forall \bar{n} \geq k \geq 0$. So we can write
\begin{equation}
	V_{\bar{n}} \leq \phi(\bar{n}-1,0) V_{0}.
\end{equation}
Under the condition of lemma 5 of \cite{guo1990estimating}, 
\begin{equation}
	\mathbf{E}[\phi(\bar{n},k)] \leq M \gamma^{n-k},
\end{equation}
$\forall \bar{n} \geq k \geq 0, 0 < \gamma < 1, M< \infty$. 
Now, for $\bar{n} \to \infty$ and $k=0$, $\gamma \to 0$ \emph{i.e.} $\phi(\bar{n}, k)) \to 0$. Hence,
\begin{equation}
	\begin{split}
		\mathbf{E}[V_{\bar{n}}] 
		\leq &~ O(\mathbf{E}[\phi(\bar{n}, 0) ||\tilde{\theta}_{0}||^{2}]) \\
		\leq &~ O(\mathbf{E}[\phi(\bar{n}, 0))] \mathbf{E}[||\tilde{\theta}_{0}||^{2}]) \to 0
	\end{split}
\end{equation}

\begin{flushright}
	$\blacksquare$	
\end{flushright}

\subsection*{B. ~ Proof of lemma 4}
%\chapter{Proof of convergence}
%\section{Proof of convergence of fixed point iteration} 
The fixed point iteration of $\hat{f}(\hat{\mathcal{X}}_{k|k})$ is derived in section 3 of \cite{chen2017maximum}. To prove the convergence of a fixed point algorithm, contraction mapping theorem which is also known as Banach fixed point theorem \cite{agarwal2001fixed} is a very important tool. Using this theorem, the convergence of proposed RMCKF can be proved.  
From section-4 of  \cite{zhao2022robust}, we can write 
\begin{equation}
	\begin{split}
		\hat{f}(\hat{\mathcal{X}}_{k|k}) = &~ (W_{k}^T \Pi_{k} W_{k})^{-1} (W_{k} \Pi_{k} D_{k}) \\
		= &~ [\hat{R}_{WW}^{G}]^{-1} [\hat{P}_{dW}^{G}].
	\end{split} \label{f_eq}
\end{equation} 
Taking the norm value, we obtain 
\begin{equation}
	||\hat{f}(\hat{\mathcal{X}}_{k|k})||_{1} \leq ||[\hat{R}_{WW}^{G}]^{-1}||_{1}  ||[\hat{P}_{dW}^{G}]||_{1}. \label{f_ineq}
\end{equation} 
Now, using the Eqn.(20, 23) of \cite{chen2017maximum}, we get
\begin{equation}
	[\hat{R}_{WW}^{G}] = \frac{1}{L} \sum_{i=1}^{L} [G_{\sigma}(e_{k,i}) w_{k,i} w_{k,i}^T],
\end{equation} and
\begin{equation}
	[\hat{P}_{dW}^{G}] = \frac{1}{L} \sum_{i=1}^{L} [G_{\sigma}(e_{k,i}) d_{k,i} w_{k,i}]. \label{P_dW}
\end{equation}
Following the Theorem 1 of \cite{chen2015convergence}, we can write
\begin{equation}
	\begin{split}
		||[\hat{R}_{WW}^{G}]^{-1}||_{1} \leq \sqrt{L} ||[\hat{R}_{WW}^{G}]^{-1}||_{2} \leq \sqrt{L} \lambda_{max} [[\hat{R}_{WW}^{G}]^{-1}]
	\end{split}, \label{r_ineq}
\end{equation}
where $\lambda_{max}[.]$ denotes the maximum eigan value of the given matrix. From (\ref{r_ineq}), we can write
\begin{equation}
	\begin{split}
		\lambda_{max} [[\hat{R}_{WW}^{G}]^{-1}] = &~ \frac{1}{\lambda_{min} [[\hat{R}_{WW}^{G}]^{-1}]} \\
		= &~ \frac{L}{\lambda_{min} [ \sum_{i=1}^{L} [G_{\sigma}(e_{k,i}) w_{k,i} w_{k,i}^T]]} \\
		\leq &~ \frac{L}{\lambda_{min} [ \sum_{i=1}^{L} [G_{\sigma}(\beta ||w_{k,i}||_{1} + |d_{k,i}|) w_{k,i} w_{k,i}^T]]},
	\end{split} \label{lamda_r}
\end{equation}
where $\lambda_{min}[.]$ denotes the minimum eigan value of the given matrix and $||e_{k,i}||_{1} = ||d_{k,i} - w_{k,i} x_{k,i}||_{1}
\leq |d_{k,i}| + \beta ||w_{k,i}||_{1}$.
Now,
\begin{equation}
	\begin{split}
		||[\hat{P}_{dW}^{G}]||_{1} = &~ ||\frac{1}{L} \sum_{i=1}^{L} [G_{\sigma}(e_{k,i}) d_{k,i} w_{k,i}]||_{1}	\\
		\leq &~ \frac{1}{L} \sum_{i=1}^{L} ||[G_{\sigma}(e_{k,i}) d_{k,i} w_{k,i}]||_{1}	\\
		\leq &~ \frac{1}{L} \sum_{i=1}^{L} |d_{k,i}| ||w_{k,i}||_{1},
	\end{split} \label{p_ineq}
\end{equation}
as $||G_{\sigma}(e_{k,i})|| \leq 1$ for any $i,k$.\\
Using (\ref{f_ineq}), (\ref{r_ineq}), (\ref{lamda_r}), and  (\ref{p_ineq}) we will get
\begin{equation}
	\begin{split}
		||\hat{f}(\hat{\mathcal{X}}_{k|k})||_{1} \leq \phi(\sigma) = 
		\frac{\sqrt{L} [\frac{1}{L} \sum_{i=1}^{L} |d_{k,i}| ||w_{k,i}||_{1}] \times L}{\lambda_{min} [ \sum_{i=1}^{L} [G_{\sigma}(\beta ||w_{k,i}||_{1} + |d_{k,i}|) w_{k,i} w_{k,i}^T]]}
	\end{split} \label{phi_sig}
\end{equation}
From (\ref{phi_sig}) it can be said that $\text{ $\displaystyle \lim_{\sigma \to 0^{+}}$} \phi(\sigma) = \infty$ and  $\text{ $\displaystyle \lim_{\sigma \to \infty}$} \phi(\sigma) = \frac{\sqrt{L} [ \sum_{i=1}^{L} |d_{k,i}| ||w_{k,i}||_{1}]} {\lambda_{min} [ \sum_{i=1}^{L} [w_{k,i} w_{k,i}^T]]} = \epsilon_{1}$. Now let us assume $\phi(\sigma^*) = \beta$. Hence, for any $\sigma \geq \sigma^*$, $\phi(\sigma) \leq \beta$ will hold, that means $||\hat{f}(\hat{\mathcal{X}}_{k|k})||_{1} \leq \beta$. 

Now again from (\ref{f_eq}), we can write
\begin{equation}
	\begin{split}
		\frac{\partial \hat{f}(\hat{\mathcal{X}}_{k|k})}{\partial \mathcal{X}_{k}} = &~ \frac{\partial}{\partial \mathcal{X}_{k}} ([\hat{R}_{WW}^{G}]^{-1} [\hat{P}_{dW}^{G}]) \\
		= &~ - [\hat{R}_{WW}^{G}]^{-1} (\frac{\partial}{\partial \mathcal{X}_{k}} \hat{R}_{WW}^{G}) [\hat{R}_{WW}^{G}]^{-1} \\
		&~ \times [\hat{P}_{dW}^{G}] + [\hat{R}_{WW}^{G}]^{-1} (\frac{\partial}{\partial \mathcal{X}_{k}} \hat{P}_{dW}^{G}) \\
		= &~ - [\hat{R}_{WW}^{G}]^{-1} (\frac{1}{L \sigma^2} \sum_{i=1}^{L} \mu_{2} e_{k,i} w_{k,i} G_{\sigma} (e_{k,i}) w_{k,i} w_{k,i}^T) \hat{f}(\mathcal{X}_{k}) \\
		&~ +  [\hat{R}_{WW}^{G}]^{-1} (\frac{1}{L \sigma^2} \sum_{i=1}^{L} \mu_{2} e_{k,i} w_{k,i} G_{\sigma} (e_{k,i}) d_{k,i} w_{k,i}) \\
		= &~ Z_{1} + Z_{2}
	\end{split} \label{z1z2}
\end{equation} 

Considering the inequalities $||\hat{f}(\hat{\mathcal{X}}_{k|k})||_{1} \leq \beta$ and $||G_{\sigma}(e_{k,i})||_{1} \leq 1$ $\forall i, k$, we can write 
\begin{equation}
	Z_{1} \leq \frac{\beta}{L \sigma^2} ||[\hat{R}_{WW}^{G}]^{-1}||_{1} \sum_{i=1}^{L} \mu_{2} (\beta ||w_{k,i}||_{1} + |d_{k,i}|) ||w_{k,i}||_{1} ||w_{k,i} w_{k,i}^T||_{1}  \label{z1_ineq}
\end{equation} and 

\begin{equation}
	Z_{2} \leq \frac{1}{L \sigma^2} ||[\hat{R}_{WW}^{G}]^{-1}||_{1} \sum_{i=1}^{L} \mu_{2} (\beta ||w_{k,i}||_{1} + |d_{k,i}|) ||w_{k,i}||_{1} |d_{k,i}| ||w_{k,i}||_{1} \label{z2_ineq}
\end{equation}
Using (\ref{r_ineq}), (\ref{lamda_r}), (\ref{z1z2}), (\ref{z1_ineq}) and (\ref{z2_ineq}), the following inequality can be written
\begin{equation}
	\begin{split}
		||\frac{\partial \hat{f}(\hat{\mathcal{X}}_{k|k})}{\partial \mathcal{X}_{k}}||_{1} \leq &~ \frac{\sqrt{L} \sum_{i=1}^{L} \mu_{2} (\beta ||w_{k,i}||_{1} + |d_{k,i}|) ||w_{k,i}||_{1}}{\lambda_{min} [\sum_{i=1}^{L} G_{\sigma} (\beta ||w_{k,i}||_{1} + |d_{k,i}|) w_{k,i} w_{k,i}^T]} \\ 
		&~ \times \frac{(\beta ||w_{k,i} w_{k,i}^T||_{1} + |d_{k,i}| ||x_{k,i}||_{1})}{\sigma^2} \\
		= &~ \psi(\sigma) \label{psi_sig}
	\end{split}
\end{equation}
From (\ref{psi_sig}) we can write  $\text{ $\displaystyle \lim_{\sigma \to 0^{+}}$} \psi(\sigma) = \infty$ and $\text{ $\displaystyle \lim_{\sigma \to \infty}$} \psi(\sigma) = 0$. Now consider $\psi(\sigma^{+}) = \alpha$ where $0 \leq \alpha \leq 1$. Therefore for any $\sigma \geq \sigma^{+}$, $\psi(\sigma) \leq \alpha$. Hence it is proved that $|| \frac{\partial \hat{f}(\hat{\mathcal{X}}_{k|k})}{\partial \mathcal{X}_{k}}||_{1} \leq \alpha \leq 1$.
\begin{flushright}
	$\blacksquare$	
\end{flushright}
%\end{appendices}

\end{document}